\documentclass[a4paper,12pt]{amsart}
\usepackage{amsmath,amssymb}
\usepackage{amssymb}
\usepackage{color}
\usepackage{enumitem}
\usepackage[nospace]{cite}

\usepackage[margin={2.8cm,3.5cm}]{geometry}
\usepackage{placeins}

\usepackage{graphicx}
\graphicspath{ {images/} }

\newcommand{\N}{\mathbb{N}}
\newcommand{\R}{\mathbb{R}}
\newcommand{\C}{\mathbb{C}}

\newcommand{\eps}{\varepsilon}

\newcommand{\ffi}{\varphi}

\newcommand{\om}{\omega}

\DeclareMathOperator \dist{dist}
\DeclareMathOperator \supp{supp}

\numberwithin{equation}{section}
 
\theoremstyle{plain}
\newtheorem{theorem}{{Theorem}}[section] 
\newtheorem{proposition}[theorem]{Proposition}
\newtheorem{lemma}[theorem]{Lemma}
\newtheorem{remark}[theorem]{Remark}
\makeatletter

\@addtoreset{equation}{section}  
\makeatother

\renewcommand{\leq}{\leqslant}	\renewcommand{\geq}{\geqslant}
\renewcommand\over[2]{{\,\buildrel #1\over#2\,}}

\newcommand{\nr}[1]{\left\Vert #1\right\Vert}         

\renewcommand{\Re}{\mathop{\rm{Re}}\nolimits}        
\renewcommand{\Im}{\mathop{\rm{Im}}\nolimits}        
\renewcommand{\a}{\alpha}\newcommand{\G}{\Gamma}\renewcommand{\d}{\delta}
 \renewcommand{\th}{\theta}\renewcommand{\l}{\lambda}\newcommand{\m}{\mu}

\renewcommand{\o}{\omega}

\usepackage{mathrsfs}

\usepackage{accents}
\makeatletter
\newcommand{\sbullet}{%
  \hbox{\fontfamily{lmr}\fontsize{.4\dimexpr(\f@size pt)}{0}\selectfont\textbullet}}

\makeatother

\begin{document}

\title[Existence]{Existence and orbital stability of standing waves to a nonlinear Schr\"odinger equation with inverse square potential on the half-line}

\author[E. Csobo]{Elek Csobo}
\address{Goethe-Universit\" at Frankfurt \\
Institut f\" ur Mathematik\\
Robert-Mayer-Str. 10\\
60629 Frankfurt am Main, Germany}
\email{Csobo@math.uni-frankfurt.de}


\begin{abstract} We investigate the properties of standing waves to a nonlinear Schr\"odinger equation with inverse-square potential on the half-line. We first establish the existence of standing waves. Then, by a variational characterization of the ground states, we establish the orbital stability of standing waves for mass sub-critical nonlinearity.  Owing to the non-compactness and the absence of translational invariance of the problem,  we apply a profile decomposition argument. We obtain convergent minimizing sequences by comparing the problem to the problem at "infinity" ({\it i.e.}, the equation without inverse square potential). Finally, we establish orbital instability by a blow-up argument for mass super-critical nonlinearity.

\end{abstract}

\maketitle

\section{Introduction}

We study the existence and orbital stability of standing waves for the following nonlinear Schr\"odinger equation with inverse square potential on the half line
\begin{equation} \label{eq-wave}
\begin{cases}\displaystyle
i u_t + u'' + c \frac{u}{x^2}+|u|^{p-1}u=0, \\
u(0)= u_0 \in H^1_0(\R^+),
\end{cases}
\end{equation}
where $u: \R\times \R^+ \rightarrow \C$, $u_0: \R^+ \rightarrow \C$, $1<p<\infty$, and $0<c<1/4$.


There has been considerable interest recently in the study of the Schr\"odinger equation with inverse-square potential in three and higher dimensions. Classification of the so-called minimal mass blow-up solutions, global well-posedness, and stability of standing wave solutions were studied in \cite{CsoboIV, BensouilahIV, DinhIV,TrachanasIV}. In the papers by A.~Bensouilah, V.~D.~Dinh, and S.~Zhu \cite{BensouilahIV}, and by G.~P.~Trachanas and N.~B.~Zographopoulos \cite{TrachanasIV} the authors establish orbital stability of ground state solutions in the Hardy subcritical $(c<(N-2)^2/4)$ and Hardy critical $(c=(N-2)^2/4)$ case respectively for dimensions higher that three. In both cases, orbital stability is proved by showing the precompactness of minimizing sequences of the energy functional on an $L^2$ constraint. Local well-posedness was established for the two-dimensional space by T.~Suzuki in \cite{SuzukiIV}, and in three and higher dimensions by N.~Okazawa, T.~Suzuki, and T.~Yokota in \cite{OkazawaIV}. The presence of the inverse square potential in one-dimensional space has also attracted attention. In \cite{KovarikIV} the H.~Kovarik and F.~Truc established dispersive estimates for $\partial_x^2+c/x^2$.

The dynamics of the equation is closely related to Hardy's inequality (see \cite{DaviesIV})
\begin{equation}\label{Hardy}
c\int_0^\infty\frac{|u|^2}{x^2}dx\leq \int_0^\infty |u'|^2dx \textit{ for all  } u \in C^\infty_0(0,\infty),
\end{equation}
where $c\leq 1/4$. We introduce the Hardy functional
\[
H(u)=\int_0^\infty\left( |u'|^2-\frac{c}{x^2}|u|^2\right)dx,
\]
which is closely related to our problem. We will mainly focus on the case $0<c<1/4$, when the natural energy space associated to \eqref{eq-wave} is $H^1_0(\R^+)$, and the semi-norm $\nr{u'}^2_{L^2}$ is equivalent to $H(u)$. 

Let us consider the operator
\[
H_c=-\frac{\partial^2}{\partial x^2} - \frac{c}{x^2}
\]
acting on $C^\infty_0(\R^+)$. Owing to the Hardy inequality, if $c<1/4$ the quadratic form $\langle H_c\ffi,\ffi\rangle$ is positive definite on $C^\infty_0(\R^+)$. It is natural to take the Friedrichs extension of $H_c$, thereby defining a self-adjoint operator in $L^2(\R^+)$, which generates an isometry group in $H^1_0(\R^+)$.

Local well-posedness for parameters $1<p<\infty$ and $0<c<\frac{1}{4}$ follows by standard arguments (see e.g. in \cite{CazenaveIV} Chapter 4). In particular, the following holds.

\begin{theorem} \label{loc-posedness}
Let $1<p<\infty$ and $c<1/4$. For any initial value $u_0 \in H^1_0(\R^+)$, there exist $T_{\mathrm{min}},T_{\mathrm{max}} \in (0,\infty]$ and a unique maximal solution $u\in C((-T_{\mathrm{min}},T_{\mathrm{max}}),H^1_0(\R^+))$ of \eqref{eq-wave}, which satisfies for all $t\in (-T_{\mathrm{min}},T_{\mathrm{max}})$ the conservation laws
\begin{equation}\label{massIV}
\nr{u(t)}_{L^2}=\nr{u_0}_{L^2}, \quad E(u(t))=E(u_0),
\end{equation}
where the energy is defined as
\begin{equation}\label{energyIV}
E(u)=\frac{1}{2}\nr{u'}^2_{L^2}-\frac{c}{2}\nr{\frac{u}{x}}^2_{L^2}-\frac{1}{p+1}\nr{u}^{p+1}_{L^{p+1}}, \textit{ for } u\in H^1_0(\R^+).
\end{equation}
Moreover, the so-called blow-up alternative holds: if ${T_\mathrm{max}}<\infty$ then $\lim_{t\rightarrow T_{\mathrm{max}}}\nr{u'(t)}_{L^2}=\infty$, (or $T_{\mathrm{min}}<\infty$ then $\lim_{t\rightarrow -T_{\mathrm{min}}}\nr{u'(t)}_{L^2}=\infty$).
\end{theorem}

In this work we address the existence of standing wave solutions and their orbital stability/instability. By introducing the ansatz $u(t,x)=e^{i\om t}\ffi(x)$, the standing wave profile equation to \eqref{eq-wave} reads as
\begin{equation} \label{st}
\ffi''-\frac{c}{x^2}\ffi+\om\ffi-|\ffi|^{p-1}\ffi=0.
\end{equation}
First we will prove regularity of standing waves and the Pohozaev identities. To establish the existence of standing waves we carry out a minimization procedure on the Nehari manifold for the so-called {\it action functional}
\[
S(v)=\frac{1}{2}\nr{v'}^2_{L^2}-\frac{c}{2}\nr{\frac{v}{x}}^2_{L^2}+\frac{\om}{2}\nr{v}^2_{L^2}-\frac{1}{p+1}\nr{v}^{p+1}_{L^{p+1}}\quad v\in H^1_0(\R^+).
\]
Owing to the non-compactness of the problem, we have to use a profile decomposition lemma, in the spirit of the article by L.~Jeanjean and K.~Tanaka \cite{JeanjeanIV}. To establish strong convergence of the minimizing sequence on the Nehari manifold we compare the minimization problem with the problem "at infinity", i.e. when $c=0$. Hence, we obtain that the set of {\it bound states} is not empty:
\[
\mathcal{A}=\{ u\in H^1_0(\R^+)\setminus \{0\} : u''+c u/x^2 -\o u + |u|^{p-1}u=0\}\neq \varnothing.
\]
We are in particular interested in the orbital stability/instability of {\it ground states}, i.e., solutions which minimize the action
\[
\mathcal{G}=\{u\in \mathcal{A} : S(u)\leq S(v) \text{ for all } v \in \mathcal{A} \}.
\] 
We use Lions' concentration-compactness principle to obtain a variational characterization of ground states on an $L^2$-constraint, thereby establishing the orbital stability of the set of ground states for nonlinearities with power $1<p<5$. Finally, for $p\geq 5$ we establish strong instability by a convexity argument. In conclusion we prove the following theorem:

\begin{theorem}
Let $0<c<1/4$, $0<\om$ and  $1<p$. Then \eqref{st} admits a nontrivial solution in $H^1_0(\R^+)$, which decays exponentially at infinity. Moreover, if $1<p<5$ the set of ground states is orbitally stable, and if $p\geq 5$, the set of ground states is strongly unstable by blow-up.
\end{theorem}

\section{Existence of bound states}

We start by investigating the standing wave equation,
\begin{equation} \label{statIV}
\begin{cases}
\ffi'' + \frac{c}{x^2}\ffi - \o\ffi + |\ffi|^{p-1}\ffi=0,\\
\ffi\in H^1_0(\R^+) \setminus \{0\}.
\end{cases}
\end{equation} 
First, we prove the regularity of solutions to \eqref{statIV} by a bootstrap argument.

\begin{proposition} \label{regularity}
Let $\o>0$ and $c < 1/4$. Assume $\ffi\in H^1_0(\R^+)$ is a solution of \eqref{statIV} in $H^{-1}(\R^+)$. Then the following statements are true
\begin{enumerate}
\item $\ffi \in W^{2,r}_0((\epsilon,\infty))$ for all $r\in [2, +\infty)$ and $\epsilon>0$, in particular $\ffi \in C^1((\epsilon,\infty))$;
\item The solution is exponentially bounded, that is $\mathrm{e}^{\sqrt{\o} x}(|\ffi|+|\ffi'|)\in L^{\infty}(\R^+)$;
\end{enumerate}
\end{proposition}

\begin{proof}
{\it (1) } For $\ffi \in H^1_0(\R^+)$ we have $\ffi \in L^q(\R^+)$ for all $q\in [2,\infty]$. We get easily that $|\ffi|^{p-1}\ffi\in L^q(\R^+)$ for all $q\in [2,\infty)$. By \eqref{statIV} we have for any $\epsilon>0 $ that $\ffi\in W^{2,q}_0((\epsilon,\infty))$ for all $q\in [2,\infty)$. 
By Sobolev's embedding we get $\ffi \in C^{1,\delta}((\epsilon,\infty))$ for all $\delta \in (0,1)$, hence $|\ffi(x)|\rightarrow 0$, and $|\ffi'(x)|\rightarrow 0$ as $x\rightarrow \infty$. 

{\it (2) } Let $\o>0$. Changing $\ffi(x)$ to $\ffi(x)=\om^{1/(p-1)}\ffi(\sqrt{\om}x) $ we may assume that $\om=1$ in \eqref{statIV}. Let $\eps>0$ and $\th_\eps(x)=e^{\frac{x}{1+\eps x}}$, for $x\geq 0$. It is easy to see that $\th_\eps$ is bounded, Lipschitz continuous, and $|\th'_\eps(x)|\leq \th_\eps(x)$ for all $x\in \R^+$. Additionally, $\th_\eps(x)\rightarrow e^x$ uniformly on bounded sets of $\R^+$. Taking the scalar product of the equation \eqref{statIV} with $\th_\eps \ffi \in H^1_0(\R^+)$, we get
\[
\Re \int_{\R^+} \ffi'\cdot (\th_\eps \bar{\ffi})'dx - c\int_{\R^+} \th_\eps \frac{|\ffi|^2}{x^2} dx + \int_{\R^+}\th_\eps|\ffi|^2 dx = \int_{\R^+}\th_\eps |\ffi|^{p+1}dx.
\]
Using the inequality $\Re(\ffi'(\th_\eps\bar{\ffi})')\geq \th_\eps|\ffi'|^2-\th_\eps|\ffi||\ffi'|$ and
\[
\int_{\R^+} \th_\eps|\ffi||\ffi'|dx \leq \frac{1}{2}\int_{\R^+}\th_\eps|\ffi|^2dx +\frac{1}{2}\int_{\R^+}\th_\eps|\ffi'|^2dx,
\]
we obtain
\[
\frac{1}{2}\int_{\R^+}\th_\eps |\ffi'|^2dx + \frac{1}{2}\int_{\R^+}\th_\eps|\ffi|^2dx - c \int_{\R^+} \th_\eps \frac{|\ffi|^2}{x^2}dx \leq \int_{\R^+} \th_\eps |\ffi|^{p+1} dx.
\]
Let $R>0$ such that if $x>R$, then $\frac{c}{x^2}\leq \frac{1}{8}$ and $|\ffi(x)|^{p-1}\leq \frac{1}{8}$. Then we get easily that
\[
c\int_{\R^+}\th_\eps \frac{|\ffi|^2}{x^2} dx + \int_{\R^+} \th_\eps |\ffi|^{p+1} \leq e^R \left(\int_0^R c\frac{|\ffi|^2}{x^2} dx +\int_0^R|\ffi|^{p+1}dx \right) + \frac{1}{4}\int_{\R^+} \th_\eps |\ffi|^2dx.
\]
From the last two inequalities it follows that
\[
\frac{1}{2}\int_{\R^+}\th_\eps |\ffi'|^2dx + \frac{1}{4}\int_{\R^+}\th_\eps|\ffi|^2dx \leq e^R \left(\int_0^R c\frac{|\ffi|^2}{x^2} dx +\int_0^R|\ffi|^{p+1}dx \right).
\]
By taking $\eps\downarrow 0$ we get
\[
\frac{1}{2}\int_{\R^+} e^{x} |\ffi'|^2dx + \frac{1}{4}\int_{\R^+}e^{x}|\ffi|^2dx <\infty.
\]
Since both $\ffi$ and $\ffi'$ are Lipschitz continuous we deduce that $|\ffi(x)|e^{x}$ and $|\ffi'(x)|e^{x}$ are bounded.
\end{proof}

We now prove that there exists a solution to \eqref{statIV}. We define the action functional associated to \eqref{statIV} as follows
\[
S(u)= \frac{1}{2} H(u) + \frac{\o}{2}\nr{u}^2_{L^2} - \frac{1}{p+1}\nr{u}^{p+1}_{L^{p+1}},
\]
for $c<1/4$ and $u\in H^1_0(\R^+)$. Clearly, we have
\[
S'(u)= - u''-\frac{c}{x^2}u + \o u - |u|^{p-1}u.
\]
Therefore, to prove the existence of a solution to \eqref{statIV} amounts to show that $S$ has a nontrivial critical point.
A simple calculation yields the following identities.
\begin{lemma} \label{identity}
Assume $p>1$, $\o>0$ and $c< 1/4$. Let $\ffi\in H^1_0(\R^+)$ be a solution of \eqref{statIV} in $H^{-1}(\R^+)$. Then the following identities are true:
\begin{align} 
\nr{\ffi'}^2_{L^2} - c\nr{\frac{\ffi}{x}}^2_{L^2} + \o\nr{\ffi}^2_{L^2}-\nr{\ffi}^{p+1}_{L^{p+1}}=0, \label{id1IV} \\
\nr{\ffi'}^2_{L^2} - c\nr{\frac{\ffi}{x}}^2_{L^2} - \frac{p-1}{2(p+1)}\nr{\ffi}^{p+1}_{L^{p+1}} = 0. \label{pohozaevIV}
\end{align}
\end{lemma}

\begin{proof}
We obtain the first equality by multiplying \eqref{statIV} by $\bar{\ffi}$ and integrating over $\R^+$.

To prove the second equality, let us put $\ffi_{\l}(x)=\l^{1/2}\ffi(\l x)$ for $\l>0$. We have that
\[
S(\ffi_\l)=\frac{\l^2}{2}\nr{\ffi'}^2_{L^2}-\frac{\l^2 c}{2}\nr{\frac{\ffi}{x}}^2_{L^2}+\frac{\o}{2}\nr{\ffi}^2_{L^2}-\frac{\l^{(p-1)/2}}{p+1}\nr{\ffi}^{p+1}_{L^{p+1}},
\]
from which we get
\[
\frac{\partial}{\partial \l} S(\ffi_\l)\Big|_{\l=1}=\nr{\ffi'}^2_{L^2} - c \nr{\frac{\ffi}{x} }_{L^2}^2-\frac{p-1}{2(p+1)}\nr{\ffi}^{p+1}_{L^{p+1}}.
\]
We also have that
\[
\frac{\partial}{\partial\l}S(\ffi_\l)\Big|_{\l=1} = \left\langle S'(\ffi), \frac{\partial\ffi_\l}{\partial \l}\Big|_{\l=1} \right\rangle.
\]
Now $\frac{\partial \ffi_\l}{\partial \l}\Big|_{\l=1}=\frac{1}{2} \ffi + x \ffi'$ is in $H^1(\R^+)$, since $\ffi$ and $\ffi'$ are exponentially decaying at infinity by Proposition \ref{regularity}. We obtain that the right hand-side is well-defined.
Since $\ffi$ is a critical point of $S$, we obtain $S'(\ffi)=0$, which concludes the proof.
\end{proof}

\begin{remark}
Since \eqref{id1IV} and \eqref{pohozaevIV} hold for solutions of \eqref{statIV}, it follows for $\om\leq 0$ that
\[
\o\nr{\ffi}^2_{L^2}=\frac{p+3}{2(p+1)}\nr{\ffi}^{p+1}_{L^{p+1}}>0.
\]
Hence, non-trivial solution of \eqref{statIV} exists only if $\o>0$.
\end{remark}

Let us define for all $u\in H^1_0(\R^+)$ the following functional:
\[
J(u)= (S'(u),u)_{H^{-1},H^1_0} =H(u)+ \o\nr{u}^2_{L^2}-\nr{u}^{p+1}_{L^{p+1}}.
\]
It follows from Lemma \ref{identity}, that $\mathcal{N}= \{ u\in H^1_0(\R^+)\setminus \{0\} : J(u)=0 \}$ contains all nontrivial critical points of $S$.
We aim to show that the infimum of the following minimization problem is attained
\begin{equation} \label{minimization1}
m=\inf \{S(u) : u \in \mathcal{N}\}=\frac{p-1}{2(p+1)}\inf\{\nr{u}^{p+1}_{L^{p+1}}:u \in \mathcal{N} \}.
\end{equation}
First we prove the following lemma.

\begin{lemma} \label{Nehari}
$\mathcal{N}$ is nonempty, and $m>0$.
\end{lemma}

\begin{proof}
Let $u\in H^1_0(\R^+)\setminus \{0\}$. Take
\[
t(u)=\left(\frac{H(u)+\o\nr{u}_{L^2}^2}{\nr{u}^{p+1}_{L^{p+1}}}\right)^{1/(p-1)}.
\]
By simple calculation, we get that $J(t(u)u)=0$, hence $t(u)u\in \mathcal{N}$. We see that
\[
m=\inf_{u\in \mathcal{N}} S(u)=\inf_{u\in \mathcal{N}}\left( S(u)-\frac{1}{p+1}J(u) \right)=\frac{p-1}{2(p+1)} \inf_{u\in \mathcal{N}} (H(u)+\o\nr{u}^2_{L^2}).
\]
It follows from Sobolev's and Hardy's inequalities, that there exists $C>0$ such that
\[
H(u) + \o\nr{u}^2_{L^2} =\nr{u}^{p+1}_{L^{p+1}}  \leq C (H(u) + \o \nr{u}^2_{L^2})^{(p+1)/2},
\]
for all $u \in \mathcal{N}$. Hence,
\[
\left(\frac{1}{C}\right)^{2/(p-1)}\leq H(u) + \o \nr{u}^2_{L^2} \text{ for all } u\in H^1_0(\R^+),
\]
which implies that
\[
m\geq \frac{p-1}{2(p+1)}\left( \frac{1}{C}\right)^{2/(p-1)}>0.
\]
\end{proof}

\begin{lemma} \label{positivity}
Let $c<1/4$, and $p>1$. Then if $u\in H^1_0(\R^+)$ is a minimizer of \eqref{minimization1}, then $|u|$ is also a minimizer. In particular, we can search for the minimizers of \eqref{minimization1} among the non-negative, real-valued functions of $H^1_0(\R^+)$.
\end{lemma}

\begin{proof}
Let $u\in H^1_0(\R^+)$ be a solution of the minimization problem \eqref{minimization1}. It is well-known that if $u\in H^1_0(\R^+)$ then $|u|\in H^1_0(\R^+)$ and $\nr{|u|'}_{L^2}\leq \nr{u'}_{L^2}$. Moreover, $\nr{|u|}_{L^{p+1}}=\nr{u}_{L^{p+1}}$. Therefore, $J(|u|)\leq J(u)$. Hence there exists a $\l\in (0,1]$ such that $J(\l |u|)=J(u)=0$. Then
\[
m\leq S(\l |u|)=\frac{p-1}{2(p+1)} \nr{\l u}^{p+1}_{L^{p+1}}\leq \frac{p-1}{2(p+1)} \nr{u}^{p+1}_{L^{p+1}}=m.
\]
Hence $\l=1$, $J(|u|)=0$, and $S(|u|)=m$.
\end{proof}

Let $m\in \R$. We say that $\{u_n\}_{n\in\N}$ is a Palais-Smale sequence for $S$ at level $m$, if
\[
S(u_n)\rightarrow m, \quad S'(u_n)\rightarrow 0 \textit{ in } H^{-1}(\R^+),
\]
as $n\rightarrow \infty$.

\begin{lemma} \label{palais-smale}
Let $c<1/4$, and $p>1$. There exists a bounded Palais-Smale sequence $\{u_n\}_{n\in \N}\subset \mathcal{N}$ for $S$ at the level $m$. Namely, there is a bounded sequence $\{u_n\}_{n\in\N}\subset \mathcal{N}$ such that, as $n\rightarrow \infty$,
\[
S(u_n)\rightarrow m, \quad S'(u_n)\rightarrow 0 \textit{ in } H^{-1}(\R^+).
\]
\end{lemma}

\begin{proof}
Since $\mathcal{N}$ is a closed manifold in $H^1_0(\R^+)$, it is a complete metric space. Hence, Ekeland's variational principle (see pp. 51-53 in \cite{StruweIV}) directly yields the existence of a Palais-Smale sequence at level $m$ in $\mathcal{N}$.  

We now show that if $\{u_n\}_{n\in\N} \subset \mathcal{N}$ and $\nr{u_n}^2_{H^1} \rightarrow \infty$, then $S(u_n)\rightarrow \infty$. Indeed, since $u_n\in\mathcal{N}$ from Hardy's inequality we get that
\[
S(u_n)=\frac{p-1}{2(p+1)}(H(u_n)+\om\nr{u_n}^2_{L^2})\geq \frac{p-1}{2(p+1)}(\min\{1,(1-4c)\}\nr{u'_n}^2_{L^2}+\om\nr{u_n}^2_{L^2}).
\] Therefore, any Palais-Smale sequence $\{u_n\}_{n\in\N}$ is bounded in $H^1_0(\R^+)$.
\end{proof}

Before proceeding to our next lemma, let us recall some classical results, see e.g. \cite{CazenaveIV}, concerning the case $c=0$. It is well-known that the set of solutions of 
\begin{equation} \label{classic}
q''-\om q+|q|^{p-1}q=0, \quad \om>0, \quad q \in H^1(\R) 
\end{equation}
is given by $\{e^{i\theta}q(\cdot + y): y\in\R, \theta \in\R  \}$, where $q$ is a symmetric, positive solution of \eqref{classic}, explicitly given by
\begin{equation}\label{explicitIV}
q(x)=\left(\frac{(p+1)\om}{2}\text{sech}^2\left(\frac{(p-1)\sqrt{\om}}{2}x \right)\right)^{1/(p-1)}.
\end{equation}
Moreover, up to translation and phase invariance, it is the unique solution of the minimization problem
\begin{align*}
m^\infty&=\inf\{S^\infty(u): u\in H^1(\R)\setminus\{0\}, J^\infty(u)=0\}\\
&=\frac{p-1}{2(p+1)}\inf\{\nr{u}^{p+1}_{L^{p+1}(\R)}: u\in H^1(\R)\setminus\{0\}, J^\infty(u)=0\},
\end{align*}
where the functionals $S^\infty$ and $J^\infty$ are defined by
\begin{align*}
S^\infty(u)&=\frac{1}{2}\nr{u'}^2_{L^2(\R)}+\frac{\om}{2}\nr{u}^2_{L^2(\R)}-\frac{1}{p+1}\nr{u}^{p+1}_{L^{p+1}(\R)},\\
J^\infty(u)&=\nr{u'}^2_{L^2(\R)}+\om\nr{u}^2_{L^2(\R)}-\nr{u}^{p+1}_{L^{p+1}(\R)}.
\end{align*}

\begin{lemma} \label{sub}
Let $0<c<1/4$, and $p>1$. Then $m<m^\infty$.
\end{lemma}

\begin{proof}
It is not hard to see that $m\leq m^\infty$, we only need to prove that $m\neq m^\infty$.
Let us first note that if $u\in H^1_0(\R^+)\setminus \{0\}$ and $J(u)<0$,  then $m<\tilde{S}(u)$, where
\[
\tilde{S}(u)=\frac{p-1}{2(p+1)}\left(H(u)+\om\nr{u}^2_{L^2} \right).
\]
Indeed, if $J(u)<0$, then let us define
\[
t(u)=\left(\frac{H(u)+\om\nr{u}^2_{L^2}}{\nr{u}^{p+1}_{L^{p+1}}}\right)^{1/(p-1)}.
\]
Hence $t(u) \in (0,1)$, $t(u) u\in \mathcal{N}$, and
\[
m\leq \tilde{S}(t(u) u)=t^2(u)\tilde{S}(u)< \tilde{S}(u).
\]
Now let us define $\psi_A(x)=q(x+A)-q(x-A)$ for $x\geq 0$. For large enough $A$ we obtain the following estimates (see Lemma \ref{approximations} in the Appendix):
\begin{align*}
&\int_0^\infty |\psi_A'|^2dx= \int_{-\infty}^\infty|q'|^2dx+O\left(\left(2A+\frac{1}{\sqrt{\o}}\right)e^{-2\sqrt{\o} A}\right),\\
& \int_0^\infty |\psi_A|^2dx= \int_{-\infty}^\infty |q|^2dx +O\left(\left(2A+\frac{1}{\sqrt{\o}}\right)e^{-2\sqrt{\o} A}\right),\\
&\int_0^\infty \frac{|\psi_A|^2}{x^2}dx \leq\frac{4}{A^2}\int_{-\infty}^\infty|q|^2dx +O\left(\frac{1}{A^2}e^{-\sqrt{\o} A}\right),\\
&\int_0^\infty |\psi_A|^{p+1}dx = \int_{-\infty}^\infty|q|^{p+1}dx + O\left( e^{-2\sqrt{\o} A}\right).
\end{align*}

Since $0<c<1/4$, we obtain for $A>0$ large enough 
\begin{align*}
J(\psi_A)&\leq\nr{q'}^2_{L^2(\R)}+\om\nr{q}^2_{L^2(\R)}-\nr{q}^{p+1}_{L^{p+1}(\R)}-\frac{4c}{A^2}\nr{q}^2_{L^2(\R)}+O\left(\frac{1}{A^2}e^{-\sqrt{\o} A}\right)\\
&=-\frac{4c}{A^2}\nr{q}^2_{L^2(\R)}+O\left(\frac{1}{A^2}e^{-\sqrt{\o} A}\right)<0,
\end{align*}
and
\begin{align*}
\tilde{S}(\psi_A)&\leq \frac{p-1}{2(p+1)}\left(\nr{q'}^2_{L^2(\R)}+\om\nr{q}^2_{L^2(\R)}-\frac{4c}{A^2}\nr{q}^2_{L^2(\R)}\right)+O\left(\frac{1}{A^2}e^{-\sqrt{\o} A}\right)\\
&=m^\infty-\frac{p-1}{2(p+1)}\frac{4c}{A^2}\nr{q}^2_{L^2(\R)}+O\left(\frac{1}{A^2}e^{-\sqrt{\o} A}\right)< m^\infty.
\end{align*}
Since $J(\psi_A)<0$, we get
\[
m<\tilde{S}(\psi_A)< m^\infty,
\]
which concludes the proof.
\end{proof}

We need the following lemma, which describes the behavior of bounded Palais-Smale sequences. We note that $H^1_0(\R^+)$ functions can be extended to functions in $H^1(\R)$ by setting $u\equiv 0$ on $\R^-$. The proof of the following statement is presented in the appendix.
\begin{lemma} \label{decomposition}
Let $\{u_n\}_{n\in\N}\subset H^1_0(\R^+)$ be a bounded Palais-Smale sequence for $S$ at level $m$. Then there exists a subsequence still denoted by $\{u_n\}_{n\in\N}$, a $u_0 \in H^1_0(\R^+)$ solution of \begin{align*}\ffi''+\frac{c}{x^2}\ffi-\o\ffi+|\ffi|^{p-1}\ffi=0,
\end{align*}
an integer $k\geq 0$, $\{x_n^i\}_{i=1}^k\subset \R^+$, and nontrivial solutions $q_i$ of \eqref{classic} satisfying
\begin{align*}
&u_n \rightharpoonup u_0 \quad  \textit{weakly in}  \quad H^1_0(\R^+), \\
&S(u_n)\rightarrow S(u_0)+\sum_{i=1}^kS^\infty(q_i), \\
&u_n-(u_0+\sum_{i=1}^kq_i(x-x_n^i))\rightarrow 0 \quad \textit{strongly in}\quad H^1(\R),\\
&|x_n^i|\rightarrow \infty,\quad |x_n^i-x_n^j|\rightarrow \infty \quad for \quad 1\leq i\neq j\leq k,
\end{align*}
where in case $k=0$, the above holds without $q_i$ and $x_n^i$.
\end{lemma}
We only need to show that the critical point of $S$ provided by Lemma \ref{decomposition} is non-trivial.

\begin{theorem} \label{existenceIV}
Let $0<c<1/4$. Then there exists $u\in \mathcal{N}\setminus \{0\}$, $u\geq 0$ a.e., such that $S(u)=m$.
\end{theorem}

\begin{proof}
We only have to prove that the $\{u_n\}_{n\in\N}$ bounded Palais-Smale sequence obtained in Lemma \ref{palais-smale} admits a strongly convergent subsequence. Assume that it is not the case. Using Lemma \ref{decomposition} we see that $k\geq 1$ and $u_n$ is weakly convergent to $u_0$ in $H^1_0(\R^+)$ up to a subsequence. Then
\begin{align*}
m=\lim_{n\rightarrow\infty} S(u_n)\geq S(u_0)+S^\infty(q)=S(u_0)+m^\infty.
\end{align*}
Now, $S(u_0)\geq 0$  since $J(u_0)=0$. Thus $m\geq m^\infty$, which contradicts Lemma \ref{sub}. Hence $k=0$ and $u_n\rightarrow u_0$ in $H^1_0(\R^+)$.
\end{proof}

\begin{lemma} \label{mass-of-GS}
Let $p>1$ and $\o>0$. There exists a $\m>0$ such that
\[
\int_0^\infty |u|^2dx= \mu, \text{ for every } u \in \mathcal{G}.
\]
The mass of ground state solutions is $\mu=\frac{m}{\o}\frac{p+3}{p-1}$. Moreover, we have
\[
\nr{u}^{p+1}_{L^{p+1}}=\frac{2(p+1)}{p-1}m, \text{ and } H(u)=m \text{ for every } u\in \mathcal{G}.
\]
\end{lemma}

\begin{proof}
Since $u \in \mathcal{G}$ is a solution of \eqref{statIV}, it satisfies \eqref{id1IV} and \eqref{pohozaevIV}. By subtracting the two identities we get
\begin{equation}\label{cs1}
\om \nr{u}_{L^2}^2=\frac{p+3}{2(p+1)}\nr{u}^{p+1}_{L^{p+1}}.
\end{equation}
Additionally, since $u$ is a ground state solution, it also solves the minimization problem \eqref{minimization1}. From \eqref{minimization1} and \eqref{pohozaevIV} we get
\begin{equation}\label{cs2}
\om \nr{u}^2_{L^2}+\frac{p-5}{2(p+1)}\nr{u}_{L^{p+1}}^{p+1}=2m.
\end{equation}

From \eqref{cs1} and \eqref{cs2} it follows

\[
\nr{u}^2_{L^2}=\frac{m}{\o}\frac{p+3}{p-1}>0.
\]
Thus, let $\mu=\frac{m}{\o}\frac{p+3}{p-1}$.
Now it follows from \eqref{minimization1} and \eqref{pohozaevIV} that 
\[
\nr{u}^{p+1}_{L^{p+1}}=\frac{2(p+1)}{p-1}m, \text{ and  } H(u)=m \text{ for every } u\in \mathcal{G}.
\]
which concludes the proof. 
\end{proof}

\section{Stability}

In this section we consider nonlinearities with $1<p<5$. Our aim is to prove orbital stability of the standing waves. To do so, we investigate the minimization problem:
\begin{equation}\label{min-erg}
I=\inf \{E(u) : u \in \Gamma\},
\end{equation}
where
\[
\Gamma = \{ u\in H^1_0(\R^+) : \nr{u}_{L^2}^2=\mu\}.
\]
and the energy $E$ is defined by \eqref{energyIV}. We will rely on a of Lions' concentration-compactness principle \cite{LionsIV} and the arguments by Cazenave and Lions \cite{Cazenave1IV}, see also in \cite{CazenaveIV}. The main problem is to obtain compactness of minimizing sequences owing to the absence of translation invariance. We define the problem at infinity by
\begin{equation}\label{min-infinity}
I^\infty=\inf\{ E^\infty(u) : u\in H^1(\R) \text{ and } \nr{u}^2_{L^2}=\mu\},
\end{equation}
where
\[
E^\infty(u)=\frac{1}{2}\int_{\R}|u'|^2dx-\frac{1}{p+1}\int_\R |u|^{p+1}dx.
\]
We recall some well-known facts about the minimization problem \eqref{min-infinity} (see \cite[Chapter 8.]{CazenaveIV}). For every $\mu >0$, there exists a unique, positive, symmetric function $q=q(\mu) \in H^1(\R)$, such that
\[
\nr{q}_{L^2}=\mu, \quad E^\infty(q)=I^\infty,
\]
and $q$ solves the nonlinear equation
\[
q''-\l q+|q|^{p-1}q=0,
\]
where $\l=\l(\mu)$.
Moreover, there exists $M>0$ such that
\[
e^{\sqrt{\l} |x|}|q(x)|\leq M \textit{ and } e^{\sqrt{\l} |x|}|q'(x)|\leq M.
\]
We proceed by proving the following lemma:

\begin{lemma} \label{strict}
If $0<c<1/4$, then the following inequality holds:
\[
I<I^\infty.
\]
\end{lemma}

\begin{proof}
For $A>0$, let $C(A)$ be a normalizing factor specified later. Let us define
\[
\Psi_A(x)=C(A)(q(x+A)-q(x-A)) \text{ for } x\geq 0.
\]
Since $q$ is even, we obtain $\Psi_A \in H^1_0(\R^+)$ and
\begin{align*}
&\int_0^\infty |\Psi_A(x)|^2dx=
C^2(A)\left(\int_{-\infty}^\infty |q|^2dx -\int_{-\infty}^\infty q(x+A)q(x-A)dx\right).
\end{align*}
We estimate the second integral by (see Lemma~\ref{approximations})
\[
\int_{-\infty}^\infty q(x+A)q(x-A)dx = O\left(\left( 2A+\frac{1}{\sqrt{\l}}\right)e^{-2\sqrt{\l} A}\right).
\]
We define
\[
C(A)=\left(\frac{\mu}{\mu-\int_{-\infty}^\infty q(x+A)q(x-A)dx}\right)^{1/2}.
\]
$C(A)$ is a continuous function of $A$, $C(A)\geq 1$, and $C(A)\rightarrow 1$ exponentially fast as $A\rightarrow \infty$. Thus, $\nr{\Psi_A}_{L^2}=\mu$ for all $A>0$.
By Lemma~\ref{approximations} in the Appendix, we obtain for $A>0$ large enough that
\begin{align*}
&\int_0^\infty |\Psi_A'|^2dx= C^2(A)\int_{-\infty}^\infty|q'|dx+O\left(\left(2A+\frac{1}{\sqrt{\l}}\right)e^{-2\sqrt{\l} A}\right),\\
&\int_0^\infty \frac{|\Psi_A|^2}{x^2}dx\leq \frac{4C^2(A)}{A^2}\int_0^\infty|\Psi_A|^2dx +O\left(\frac{1}{A^2}e^{-\sqrt{\l} A}\right),\\
&\int_0^\infty |\Psi_A|^{p+1}dx = C^{p+1}(A)\int_{-\infty}^\infty|q|^{p+1}dx + O(e^{-2\sqrt{\l} A}).
\end{align*}
Hence for $A$ large enough we get
\begin{align*}
E(\Psi_A)&=\frac{1}{2}\int_0^\infty|\Psi_A'|^2dx -\frac{c}{2}\int_0^\infty\frac{|\Psi_A|^2}{x^2}dx-\frac{1}{p+1}\int_0^\infty|\Psi_A|^{p+1}dx\\
&\leq C^2(A)\left(\frac{1}{2}\int_{-\infty}^\infty|q'|^2dx-\frac{C^{p-1}(A)}{p+1}\int_{-\infty}^\infty|q|^{p+1}dx \right)-\frac{c}{2}\frac{4C^2(A)}{A^2}\int_0^\infty|\Psi_A|^2dx\\
&+O\left(\frac{1}{A^2}e^{-\sqrt{\l} A}\right).
\end{align*}
Owing to the exponential decay of the last term, for large $A$ we get
\[
E(\Psi_A)\leq E(q) -\frac{2c}{A^2}\mu=I^\infty-\frac{2c}{A^2}\mu.
\]
Since $0<c<1/4$ we get that $E(\Psi_A)<I^\infty$, which concludes the proof.
\end{proof}

We need the following version of the concentration-compactness principle. The proof follows the same way as in the classical case (see \cite{LionsIV}).

\begin{lemma}\label{cc}
Let $0<c<1/4$, and $\{ u_n\}_{n\in\N} \subset H^1_0(\R^+)$ be a sequence satisfying
\[
\lim_{n\rightarrow \infty} \nr{u_n}^2_{L^2} = M \text{ and } \lim_{n\rightarrow\infty} H(u_n) < \infty.
\]
Then there exists a subsequence $\{u_n\}_{n\in\N}$ such that it satisfies one of the following alternatives.

(Vanishing) $\lim_{n\rightarrow \infty}\nr{u_n}_{L^p}\rightarrow 0$ for all $p\in(2,\infty)$.

(Dichotomy) There are sequences $\{v_n\}_{n\in\N}, \{w_n\}_{n\in\N}$ in $H^1_0(\R^+)$ and a constant $\a \in (0,1)$ such that:
\begin{enumerate}
	\item $\dist(\supp(v_n),\supp(w_n))\rightarrow \infty$;
	\item $|v_n|+|w_n|\leq |u_n|$;
	\item $\sup_{n\in\N}(\nr{v_n}_{H^1}+\nr{w_n}_{H^1})<\infty$;
	\item $\nr{v_n}_{L^2}^2\rightarrow \a M$ and $\nr{w_n}_{L^2}^2\rightarrow (1-\a)M$ as $n\rightarrow \infty$;
	\item $\lim_{n\rightarrow \infty}\left|\int_0^\infty |u_n|^qdx-\int_0^\infty|v_n|^qdx-\int_0^\infty|w_n|^qdx\right|=0$ for all $q\in [2,\infty)$;
	\item $\liminf_{n\rightarrow \infty}\{H(u_n)-H(v_n)-H(w_n)\}\geq 0$.
\end{enumerate}

(Compactness) There exists a sequence $y_n \in \R^+$, such that for any $\eps >0$ there is an $R>0$ with the property that
\[
\int_{(y_n-R,y_n+R)\cap \R^+} |u_n|^2\geq M-\eps.
\]
for all $n\in \N$.
\end{lemma}

We are now in a position to prove the following lemma.

\begin{lemma}\label{compactIV}
Let $1<p<5$, $0<c<1/4$, and $\om>0$. Then the infimum in \eqref{min-erg} is attained.
Additionally, all minimizing sequences are relatively compact, that is if $\{u_n\}_{n\in\N}$ satisfies $\nr{u_n}_{L^2}^2\rightarrow \mu$ and $E(u_n)\rightarrow I$ then there exists a subsequence $\{u_n\}_{n\in\N}$ which converges to a minimizer $u\in H^1_0(\R^+)$.
\end{lemma}

\begin{proof}
Step 1. We first show that $0>I>-\infty$. Let $u\in \Gamma$. For $\l>0$, we define $u_\l(x)=\l^{1/2}u(\l x) \in \Gamma$. Clearly,
\[
E(u_\l)=\frac{\l^2}{2}\nr{u'}^2_{L^2}-\frac{c\l^2}{2}\int_0^\infty\frac{|u|^2}{x^2}dx-\frac{\l^{(p-1)/2}}{p+1}\nr{u}^{p+1}_{L^{p+1}}
\]
Since $1<p<5$, we can choose a small $\l>0$ such that $E(u_\l)<0$. Hence $I<0$.

Since $c\in (0,1/4)$, we have $H(u)\sim \nr{u'}^2_{L^2}$. We get from the Gagliardo-Nirenberg inequality that there exists $C>0$ such that for all $u\in H^1_0(\R^+)$
\[
\int_0^\infty |u|^{p+1}dx\leq C H(u)^{\frac{p-1}{4}}\left(\int_0^\infty |u|^2dx\right)^{1+\frac{p-1}{4}}.
\]
Since $1<p<5$, this yields that there exists $\d>0$ and $K>0$ such that
\begin{equation} \label{lower-bound}
E(u)\geq \d \nr{u}^2_{H^1}-K \text{ for all } u \in \Gamma,
\end{equation}
from which follows that $I>-\infty$.

Every minimizing sequence is bounded in $H^1_0(\R^+)$ and bounded from below in $L^{p+1}(\R^+)$. Indeed, let $\{u_n\}_{n\in \N}\subset\Gamma$ be a minimizing sequence, then by \eqref{lower-bound} it is bounded in $H^1_0(\R^+)$. Furthermore, for $n$ large enough we have $E(u_n)<I/2$, thus
\begin{equation} \label{non-vanishing}
\nr{u_n}_{L^{p+1}}^{p+1}>-\frac{p+1}{2}I.
\end{equation}
Now $I<0$, hence the result follows.

Step 2. We now verify that all minimizing sequences have a subsequence which converges to a limit $u$ in $H^1_0(\R^+)$. Let $\{u_n\}_{n\in\N}$ satisfying $\nr{u_n}^2_{L^2}\rightarrow \mu$ and $E(u_n)\rightarrow I$. Since every minimizing sequence is bounded in $H^1_0(\R^+)$, $\{u_n\}_{n\in\N}$ has a weak-limit $u\in L^p(\R^+)$ . We can apply the concentration-compactness principle (see Lemma \ref{cc}) to the sequence $\{u_n\}_{n\in\N}$. We note that since the sequence is bounded from below in $L^{p+1}(\R^+)$ vanishing cannot occur.

Now let us assume that dichotomy occurs. Let $\a \in (0,1)$, $\{v_n\}_{n\in\N}$ and $\{w_n\}_{n\in\N}$ sequences as in Lemma \ref{cc}. It follows from {\it (5)} and {\it (6)} of Lemma \ref{cc} that
\[
\liminf_{n\rightarrow \infty}(E(u_n)-E(v_n)-E(w_n))\geq 0,
\]
hence
\begin{equation} \label{contradiction}
\limsup_{n\rightarrow \infty}(E(v_n)+E(w_n))\leq I.
\end{equation}
Observe that for $u\in H^1_0(\R^+)$, and $a>0$, we have
\[
E(u)=\frac{1}{a^2}E(au)+\frac{a^{p-1}-1}{p+1}\int_0^\infty |u|^{p+1}dx.
\]
Let $a_n=\sqrt{\mu}/\nr{v_n}_{L^2}$ and $b_k^2=\sqrt{\mu}/\nr{w_n}_{L^2}$. Hence, $a_nv_n \in \Gamma$ and $b_nw_n\in \Gamma$, which implies
\begin{align*}
E(v_n)\geq\frac{I}{a^2_n}+\frac{a_n^{p-1}-1}{p+1}\int_0^\infty|v_n|^{p+1}dx,
\\
E(w_n)\geq\frac{I}{b^2_n}+\frac{b_n^{p-1}-1}{p+1}\int_0^\infty|w_n|^{p+1}dx.
\end{align*}
Therefore
\[
E(v_n)+E(w_n)\geq I(a_n^{-2}+b_n^{-2})+\frac{a_n^{p-1}}{p+1}\int_0^\infty|v_n|^{p+1}+\frac{b_n^{p-1}}{p+1}\int_0^\infty|w_n|^{p+1}.
\]
Now we observe $a_n^{-2}\rightarrow \a$ and $b_n^{-2}\rightarrow (1-\a)$ by {\it (4)} of Lemma \ref{cc}. Since $\alpha \in (0,1)$, we get that $\theta = \min \{\alpha^{-(p-1)/2};(1-\alpha)^{-(p-1)/2}) \}>1$. Property {\it (5)} of Lemma \ref{cc} and \eqref{non-vanishing} implies
\begin{align*}
\liminf_{n\rightarrow \infty}(E(v_n)+E(w_n))\geq I+ \frac{\theta-1}{p+1}\liminf_{n\rightarrow \infty}\int_0^\infty|u_n|^{p+1}dx,
\geq I+\frac{\theta-1}{2}>I,
\end{align*}
which contradicts with \eqref{contradiction}. Hence the following holds: there exists a sequence $y_n\in\R^+$, such that for any $\eps>0$ there exists $R>0$ with the property that
\begin{equation}\label{cmp}
\int_{(y_n-R,y_n+R)\cap\R^+}|u_n|^2_{L^2} \geq \mu-\eps.
\end{equation}
for all $k\in\N$.

We now show that $\{y_n\}_{n\in\N}$ is bounded in $\R^+$. First we show that if $y_n\rightarrow \infty$, then
\begin{equation} \label{vanish}
\lim_{n\rightarrow\infty}\int_0^\infty \frac{|u_n|^2}{x^2}dx=0.
\end{equation}
Let us assume by contradiction that 
\begin{equation}\label{assumption}
\int_0^\infty \frac{|u_n|^2}{x^2}dx\geq \delta >0,
\end{equation}
which implies together with Hardy's inequality that
\begin{equation} \label{assumption1}
H(u_n)\geq (1/4-c)\delta.
\end{equation}
Let us take $\xi\in C^\infty(\R^+)$, such that for $\tilde{R}>0$ and $a>0$ we have that $\xi(r)=1$ for $0\leq r \leq \tilde{R}$, $\xi(r)=0$ for $r\geq \tilde{R}+a$, and $\nr{\xi'}_{L^\infty}\leq 2/a$. We introduce $u_{n,1}=u_n\cdot\xi$ and $u_{n,2}=u_n\cdot(1-\xi)$. Clearly, $u_{n,1}\in H^1_0(\R^+)$, $u_{n,2}\in H^1_0(\R^+)$ and $u_n=u_{n,1}+u_{n,2}$. Moreover, the following inequalities hold
\begin{align*}
|u'_{n,1}|^2\leq 2(4a^{-2}|u_n|^2+|u'_n|^2),
\\
|u'_{n,2}|^2\leq 2(4a^{-2}|u_n|^2+|u'_n|^2).
\end{align*}
We obtain by direct calculation that
\[
E(u_n)=E(u_{n,1})+E(u_{n,2})+\rho_n
\]
where
\begin{align*}
\rho_n&=\frac{1}{2}\int_{\tilde{R}}^{\tilde{R}+a} \left[(|u_n'|^2-|u'_{n,1}|^2-|u'_{n,2}|^2)-\frac{c}{x^2}(|u_n|^2-|u_{n,1}|^2-|u_{n,2}|^2) \right]dx
\\&-\frac{1}{p+1}\int_{\tilde{R}}^{\tilde{R}+a}(|u_n|^{p+1}-|u_{n,1}|^{p+1}-|u_{n,2}|^{p+1})dx.
\end{align*}
We show that there exists $\tilde{R}>0$ and $a>1$, such that for $n$ large enough $|\rho_n|\leq (1/4-c)\frac{\delta}{4}$. First we observe by the properties of the cut-off that
\[
\left|\frac{1}{2}\int_{\tilde{R}}^{\tilde{R}+a} (|u_n'|^2-|u'_{n,1}|^2-|u'_{n,2}|^2)dx\right|\leq \frac{5}{2} \int_{\tilde{R}}^{\tilde{R}+a}|u_n'|^2dx+\frac{8}{a^2}\int_{\tilde{R}}^{\tilde{R}+a}|u_n|^2dx.
\]
We claim that there exist $\tilde{R}>0$ and $a>1$ such that for a subsequence $\{u_{n_k}\}$ we have
\begin{equation}\label{boundedness1}
\int_{\tilde{R}}^{\tilde{R}+a}|u'_{n_k}|^2dx< \frac{1}{20}(1/4-c)\delta.
\end{equation}
Suppose that this claim does not hold, that is for all $R>0$, $a>1$ there exists $k\in \N$ such that for all $n\geq k$ the following holds
\[
\int_R^{R+a}|u'_n|^2dx\geq \frac{1}{20}(1/4-c)\delta.
\]
Let $(R_1,R_1+a_1)$. There exists $k_1\in \N$, such that for all $n\geq k_1$ we have
\[
\int_{R_1}^{R_1+a_1}|u'_n|^2dx\geq \frac{1}{20}(1/4-c)\delta.
\]
Now let $R_2>R_1+a_1$ and $a_2>1$. Then by our assumption there exists $k_2\in \N$, such that for all $n\geq k_2$ it holds that
\[
\int_{R_2}^{R_2+a_2}|u'_n|^2dx\geq \frac{1}{20}(1/4-c)\delta.
\]
Hence, there exists a subsequence $\{v_{n_k}\}_{k\in\N}$ such that for all $j\in\{1,2\}$ it holds that
\[
\int_{R_j}^{R_j+a_j}|u'_{n_k}|^2dx\geq \frac{1}{20}(1/4-c)\delta
\]
for all $k\in\N$. Therefore, we can construct for all $l\in\N$ a subsequence $\{u_{n_k}\}_{k\in\N}$, such that for all $1\leq j\leq l$ there are disjoint intervals $A_j=(R_j,R_j+a_j)$, such that
\[
\int_{A_j}|u'_{n_k}|^2dx\geq \frac{1}{20}(1/4-c)\delta.
\]
Hence for all $l\in\N$ there exists a subsequence $\{ u_{n_k}\}_{k\in\N}$, such that for all $k\in\N$ we have
\[
\int_0^\infty|u'_{n_k}|^2dx\geq \sum_{j=1}^{l}\int_{A_j}|u'_{n_k}|^2dx\geq \frac{l}{20}(1/4-c)\delta.
\]
This implies that $\int_0^\infty|u'_{n_k}|^2dx\rightarrow \infty$, which is a contradiction since $\{u_n\}_{n\in\N}$ is bounded in $H^1_0(\R^+)$. Hence the assertion \eqref{boundedness1} is true. Now we note that
\[
\int_0^R |u_n|^{p+1}dx \leq \nr{u_n}^{p-1}_{L^\infty}\int_0^R|u_n|^2dx.
\]
Since $\{u_n\}_{n\in\N}$ is bounded in $L^\infty(\R^+)$, in view of \eqref{cmp} we obtain for $R>0$ given in \eqref{cmp} that
\begin{equation}\label{interpolation}
\int_0^R|u_n|^2dx \rightarrow 0 \quad \text{implies} \quad \int_0^R |u_n|^{p+1}dx\rightarrow 0.
\end{equation}
For large $n$ we have $\tilde{R}+a<y_n-R$, since $y_n\rightarrow \infty$ by our assumption. Now \eqref{interpolation} implies
\begin{align}
&\left|\frac{8}{a^2}\int_{\tilde{R}}^{\tilde{R}+a}|u_n|^2dx\right|+\left|\int_{\tilde{R}}^{\tilde{R}+a}\frac{c}{x^2}(|u_n|^2-|u_{n,1}|^2-|u_{n,2}|^2) dx\right| \nonumber
\\
&+\left|\frac{1}{p+1}\int_{\tilde{R}}^{\tilde{R}+a}(|u_n|^{p+1}-|u_{n,1}|^{p+1}-|u_{n,2}|^{p+1})dx\right| \nonumber
\\
&\leq \left|\frac{8}{a^2}\int_{\tilde{R}}^{\tilde{R}+a}|u_n|^2dx\right|+\frac{c}{\tilde{R}^2}\left|\int_{\tilde{R}}^{\tilde{R}+a}|u_n|^2(1-\xi^2-(1-\xi)^2) dx\right|  \nonumber
\\
&+\left|\frac{1}{p+1}\int_{\tilde{R}}^{\tilde{R}+a}|u_n|^{p+1}(1-\xi^{p+1}-(1-\xi)^{p+1})dx\right| \nonumber
\\
&\leq \frac{(1/4-c)\delta}{8}.\label{boundedness2}
\end{align}
for large $n$. Now \eqref{boundedness1} and \eqref{boundedness2} implies 
\begin{equation} \label{estimate}
|\rho_n|\leq \frac{(1/4-c)\delta}{4}. 
\end{equation}
Let us observe that $\nr{u_{n,1}}_{L^{p+1}}\rightarrow 0$ by \eqref{interpolation}. Hence
\[
E(u_{n,1})=\frac{1}{2}H(u_{n,1})+o(1).
\]
Now let us notice that $\mathrm{supp} (u_{n,2}) \subset(\tilde{R}, \infty)$. Moreover, in view of \eqref{cmp},
\[
\int_0^\infty|u_{n,2}|^2dx = \int_{y_n-R}^\infty |u_{n,2}|^2dx+o(1).
\]
Hence
\[
\int_0^\infty \frac{|u_{n,2}|^2}{x^2}dx= \int_{y_n-R}^\infty \frac{|u_{n,2}|^2}{x^2}dx + o(1)\leq \frac{\mu}{|y_n-R|^2}.
\]
Now $y_n\rightarrow \infty$ implies that
\[
 E(u_{n,2})= E^\infty(u_{n,2})+o(1).
\]
Thus,
\[
E(u_n)=\frac{1}{2}H(u_{n,1})+E^\infty(u_{n,2})+\rho_n+o(1).
\]
From the properties of the cut-off and \eqref{cmp}, we get
\[
\nr{u_{n,2}}_{L^2}^2=\nr{u_n}^2_{L^2}-\nr{u_{n,1}}^2_{L^2}-2\Re\int_{R'}^{R'+a} u_{n,1}\bar{u}_{n,2}dx \rightarrow \mu.
\]
Since $\frac{1}{2}H(u_{n,1})+\rho_n>0$ by \eqref{assumption1} and \eqref{estimate}, we obtain
\[
I=\lim_{n\rightarrow\infty}E(u_n)\geq \lim_{n\rightarrow\infty}E^\infty(u_{n,2})\geq I^\infty.
\]
which is a contradiction, hence \eqref{vanish} follows.

Now, from \eqref{vanish} we obtain
\begin{align*}
\lim_{n\rightarrow\infty}&\left( \frac{1}{2}\int_0^\infty|u'_n|^2dx-\frac{c}{2}\int_0^\infty\frac{|u_n|^2}{x^2}dx-\frac{1}{p+1}\int_0^\infty|u_n|^{p+1}dx\right)=
\\
&=\lim_{n\rightarrow\infty} \left(\frac{1}{2}\int_0^\infty|u'_n|^2dx-\frac{1}{p+1}\int_0^\infty|u_n|^{p+1}dx\right).
\end{align*}
Hence
\[
I\geq I^\infty,
\]
which is again a contradiction. Thus $\{y_n\}_{n\in\N}$ is bounded and has an accumulation point $y^*\in \R^+$. Therefore, it follows that for any $\eps>0$ there is $R>0$ such that
\[
\int_0^R|u_n|^2_{L^2} \geq \mu-\eps.
\]
for all $k\in\N$. Hence $u_n\rightarrow u$ strongly in $L^{2}(\R^+)$. Moreover, since $\{u_n\}$ is bounded in $H^1_0(\R^+)$ it is also strongly convergent in $L^{p+1}(\R^+)$. By the weak-lower semicontinuity of $H$ (see \cite{MontefuscoIV}), it follows that $E(u)\leq \lim_{n\rightarrow \infty}E(u_n)=I$. Hence $E(u)=I$, and $E(u_n)\rightarrow E(u)$ implies that $H(u_n)\rightarrow H(u)$, which concludes that proof. 
\end{proof}


\begin{remark}
If $c<0$, the infimum is not attained on the $L^2$ constraint. Indeed, let us assume that there exists $v\in H^1_0(\R^+)$, such that $\nr{v}^2_{L^2}=\mu$ and $E(v)=I$. Then taking translates of $v$, i.e. $v(\cdot -y)$ for $y>0$, we get $E(v(\cdot-y))<I$, which is a contradiction.
\end{remark}

\begin{lemma}\label{minimization-mass}
Let $0<c<1/4$, $\o>0$ and $1<p<5$. Let $\mu$ be defined by Lemma \ref{mass-of-GS}.
Then $u\in H^1_0(\R^+)$ is a ground state solution of \eqref{statIV} if and only if $u$ solves the minimization problem
\begin{equation} \label{minmass}
\begin{cases}
u \in \Gamma,
\\
S(u)=\inf \{S(v): v\in \Gamma\}.
\end{cases}
\end{equation}
\end{lemma}

\begin{proof}
Step 1. Let us first define
\[
m_{\mathcal{A}}=\inf \{S(u) : u \in \mathcal{A}\},
\]
and
\[
m_{\Gamma}=\inf \{S(u) : u \in \Gamma\}.
\]
If $u\in\mathcal{G}$, then $S(u)=m_\Gamma$. By Lemma \ref{mass-of-GS} we know that $u\in\Gamma$, hence $m_{\mathcal{A}}\leq m_\Gamma $.
\\
Step 2. We claim that every solution of \eqref{minmass} belongs to $\mathcal{A}$. 
Indeed, let us consider a solution $u$ to \eqref{minmass}. There exists a Lagrange multiplier $\l_1\in\R$ such that $S'(u)=\l_1 u$. Hence there exists $\l\in\R$ such that
\begin{equation} \label{delta1}
-u'' -\frac{c}{x^2}u+ \l\o u= |u|^{p-1}u.
\end{equation}
Indeed, since $u$ is a solution of \eqref{minmass}, and for $\l>0$ let
\[
u_\l(x)=\l^{1/2}u(\l x).
\]
We have $u_\l \in \G$. Since $u_1$ is a solution of \eqref{minmass}, we get from \eqref{delta1} and Lemma \ref{identity} that
\begin{equation} \label{delta}
\frac{\partial}{\partial \lambda}S(u_\l)|_{\l=1}=\nr{u'}^2_{L^2}-c\nr{\frac{u}{x}}^2_{L^2}-\frac{p-1}{2(p+1)} \nr{u}^{p+1}_{L^{p+1}}=0.
\end{equation} 
We can deduce directly from \eqref{delta1} and \eqref{delta} that
\[
\l\o\m=\frac{p+3}{p-1} H(u),
\]
which implies that $\l>0$. Let us define $v$ by
\[
u(x)=\l^{1/(p-1)}v(\l^{1/2}x).
\]
By \eqref{delta}, $v\in \mathcal{A}$, hence
\[
S(v)\geq m_{\mathcal{A}}.
\]
We obtain simple calculation that
\[
m_{\Gamma}=S(u)=\l^{2/(p-1)+1/2}S(v) + (1-\l)\frac{\o\m}{2}.
\]
Hence,
\[
m_{\mathcal{A}}\geq \l^{\frac{2}{p-1}+\frac{1}{2}}m_{\mathcal{A}}+(1-\l)\frac{\o\m}{2}.
\]
Since $u$ is a solution of \eqref{delta1}, we obtain from Lemma \ref{identity} that $m_{\mathcal{A}}\geq 0$. By Lemma \ref{identity} and Lemma \ref{mass-of-GS} we have that
\[
\frac{\o\m}{2}=\left(\frac{2}{p-1}+\frac{1}{2} \right)m_{\mathcal{A}},
\]
hence
\[
1\geq \l^{\frac{2}{p-1}+\frac{1}{2}} - \l\left(\frac{2}{p-1}+\frac{1}{2} \right)+ \left(\frac{2}{p-1}-\frac{1}{2} \right).
\]
The right hand side is always strictly positive, except if $\l=1$. Thus, $\l=1$, which implies together with \eqref{delta} that $u\in \mathcal{A}$. 
\\
\\
Step 3. It follows from Step 2, that $l\leq m_{\mathcal{A}}$, hence $m_\Gamma=m_{\mathcal{A}}$. In particular, it follows that if $u\in \mathcal{G}$, then $u \in \Gamma$ and $S(u)=m_{\mathcal{A}}$, thus $u$ satisfies \eqref{minmass}. Conversely, let $u$ be the solution of \eqref{minmass}. Then by Step 2 $u\in \mathcal{A}$, and $S(u)=m_\Gamma=m_{\mathcal{A}}$, hence $u \in \mathcal{G}$.
\end{proof}

\begin{theorem}
Let $0<c<1/4$, $\om>0$, and $1<p<5$. If $\ffi$ is a ground state solution of \eqref{statIV}, then the standing wave $u(t,x)=e^{i\om t}\ffi(x)$ is an orbitally stable solution of \eqref{eq-wave}, i.e. for all $\eps>0$ there is $\delta>0$, such that if $u(0)\in H^1_0(\R^+)$ satisfies $\nr{\ffi-u(0)}_{H^1}<\delta$, then the corresponding maximal solution $u$ of \eqref{eq-wave} satisfies
\[
\sup_{t\in\R}\inf_{\theta\in\R}\nr{u(t)-e^{i\theta}\ffi}_{H^1}<\eps.
\]
\end{theorem}

\begin{proof}
Assume by contradiction that there exist a sequence $\{\ffi_n\}_{n\in \N}\subset H^1_0(\R^+)$, a sequence $\{t_n\}_{n\in\N}\subset \R$, and $\eps>0$, such that
\[
\lim_{n\rightarrow\infty} \nr{\ffi_n-\ffi}_{H^1}=0,
\]
and the corresponding maximal solution $u_n$ of \eqref{eq-wave} with initial value $\ffi_n$ satisfies
\[ 
\inf_{\theta\in\R}\nr{u_n(t_n)-e^{i\theta}\ffi}_{H^1}\geq \eps.
\]
Set $v_n=u_n(t_n)$. Applying Lemma \ref{minimization-mass}, we obtain
\begin{equation} \label{contradiction1}
\inf_{\ffi\in\mathcal{G}}\nr{v_n-\ffi}_{H^1}\geq\eps.
\end{equation}
By the conservation of charge and energy, we obtain
\[
\nr{v_n}^2_{L^2}\rightarrow \m, \textit{ and } E(v_n)\rightarrow I.
\]
Hence $\{v_n\}_{n\in\N}$ is a minimizing sequence of \eqref{min-erg}. It follows from Lemma \ref{compactIV}, that there exists a solution $u$ of the problem \eqref{min-erg}, such that $\nr{v_n-u}_{H^1}\rightarrow 0$. By Lemma \ref{minimization-mass} we obtain that $u\in\mathcal{G}$, which contradicts \eqref{contradiction1}.
\end{proof}

\section{Instability}

In this section we assume that $p \geq 5$. Let us define for $v\in H^1_0(\R^+)$ the functional
\[
Q(v)= \nr{v'}^2_{L^2}-c\nr{\frac{v}{x}}^2_{L^2}-\frac{p-1}{2(p+1)}\nr{v}^{p+1}_{L^{p+1}}. 
\]
In Lemma \ref{identity} we have shown that if $v$ is a solution of \eqref{statIV}, then $Q(v)=0$.
First, we prove the virial identities.

\begin{proposition} \label{virial}
Let $u_0\in H^1_0(\R^+)$ be such that $x u_0 \in L^2(\R^+)$ and $u$ be the corresponding maximal solution to \eqref{eq-wave}. Then  $x u(t) \in L^2 (\R^+)$ for any $t\in (-T_{\mathrm{min}},T_{\mathrm{max}})$. Moreover, the following identities hold for all $v\in H^1_0(\R^+)$:
\begin{align*}
\frac{\partial}{\partial t} \nr{x u(t)}^2_{L^2} &= 4\Im \int_0^\infty \bar{u}(t) x u'(t) dx,
\\
\frac{\partial^2}{\partial t^2} \nr{x u(t)}^2_{L^2} &= 8 Q(u(t)).
\end{align*}
\end{proposition}

\begin{proof}
The proof follows the same line as in \cite{CsoboIV}.
\end{proof}

\begin{proposition} \label{simple-blow-up}
Let $p\geq 5$ and let $u_0\in H^1_0(\R^+)$ be such that
\[
xu_0\in L^2(\R^+) \text{ and } E(u_0)<0.
\]
Then the maximal solution $u$ to \eqref{eq-wave} with initial condition $u_0$ blows up in finite time.
\end{proposition}

\begin{proof}
First, let us note that
\[
Q(u(t))=2E(u(t))+\frac{5-p}{2(p+1)}\nr{u(t)}^{p+1}_{L^{p+1}}.
\]
Since $p\geq 5$, we get by the conservation of the energy that
\[
Q(u(t))\leq 2E(u_0)<0 \text{ for all } t\in (-T_{\mathrm{min}},T_{\mathrm{max}}).
\]
Hence, Proposition \ref{virial} implies that
\[
\frac{\partial^2}{\partial t^2} \nr{xu(t)}^2_{L^2}\leq 16 E(u_0) \text{ for all } t\in (-T_{\mathrm{min}},T_{\mathrm{max}}).
\]
Integrating twice, we get
\begin{equation}\label{simple1}
\nr{xu(t)}^2_{L^2}\leq 8E(u_0)t^2 +\left(4\Im \int_0^\infty \bar{u}_0 x u'_0dx \right)t +\nr{xu_0}^2_{L^2}
\end{equation}
The main coefficient of the second order polynomial on the right hand side is negative. Thus, it is negative for $|t|$ large, what contradicts with $\nr{xu(t)}^2_{L^2}\geq 0$ for all $t$. Therefore, $-T_{\mathrm{min}}>-\infty$ and $T_{\mathrm{max}}<+\infty$.
\end{proof}

\begin{theorem} \label{virialblowup}
Assume that $\om>0$ and $p=5$. Then for any solution $\ffi\in H^1_0(\R^+)$ of \eqref{statIV} the standing wave $e^{i\o t} \ffi (x)$ is unstable by blow-up.
\end{theorem}

\begin{proof}
Since $p=5$, we have for all $v\in H^1_0(\R^+)$, that $2E(v)=Q(v)$. Hence from Lemma \ref{identity} we get that
\[
E(\ffi)=0.
\]
Let us define $\ffi_{n,0}=\left(1+\frac{1}{n}\right)\ffi$. It is easy to see that $E(\ffi_{n,0})<0$. By Lemma \ref{regularity} we know that $x \ffi_{n,0}\in L^2(\R^+)$. The conclusion follows from Proposition \ref{simple-blow-up}.
\end{proof}

\begin{theorem} \label{unstable}
Let $p>5$. Then for any ground state solution $\ffi$ to \eqref{statIV}, the corresponding standing wave $e^{i\o t}\ffi(x)$ is orbitally unstable. 
\end{theorem}

We need to prove a series of Lemmas to establish Theorem \ref{unstable}.

\begin{lemma} \label{scaling}
Let $v\in H^1_0(\R^+)\setminus \{0\}$ such that $Q(v)\leq 0$, and set $v_\l(x)=\l^{1/2}v(\l x)$ for $\l>0$. Then there exists $\l^* \in (0,1]$ such that the following assertions hold:
\begin{enumerate}
\item $Q(v_{\l^*})=0$.
\item $\l^*=1$ if and only if $Q(v)=0$.
\item $\frac{\partial}{\partial \l} S(v_\l) =\frac{1}{\l} Q(v_\l)$.
\item $\frac{\partial}{\partial \l} S(v_\l) >0$ for all $\l\in(0,\l^*)$, and $\frac{\partial}{\partial \l} S(v_\l) <0 $ for all $\l\in (\l^*, +\infty)$.
\item The function $(\l^*, +\infty)\ni\l \mapsto S(v_\l)$ is concave.
\end{enumerate}
\end{lemma}

\begin{proof}
We get that by the scaling properties of $\l \mapsto Q(v_\l)$ that
\[
Q(v_\l)=\l^2\nr{v'}^2_{L^2}-\l^2 c \nr{\frac{v}{x}}^2_{L^2}-\l^{\frac{p-1}{2}}\frac{p-1}{2(p+1)}\nr{v}^{p+1}_{L^{p+1}}.
\]
We get from the Hardy inequality that for $c\in (0,1/4)$
\[
(1-4c)\l^2\nr{v'}^2_{L^2}-\l^{\frac{p-1}{2}}\frac{p-1}{2(p+1)}\nr{v}^{p+1}_{L^{p+1}} \leq Q(v_\l)\leq\l^2\nr{v'}^2_{L^2}-\l^{\frac{p-1}{2}}\frac{p-1}{2(p+1)}\nr{v}^{p+1}_{L^{p+1}}.
\]
Since $p>5$, there exists $\l\in(0,1]$ small enough, such that $Q(v_\l)>0$. Hence, there exists $\l^*\in(0,1]$, such that $Q(v_{\l^*})=0$. This proves {\it (1)}. To prove {\it (2)}, we first note that if $\l_0=1$, then clearly $Q(v)=0$. Now assume that $Q(v)=0$. Then
\begin{align*}
Q(v_\l)&=\l^2Q(v)+(\l^2-\l^{\frac{p-1}{2}})\frac{p-1}{2(p+1)}\nr{v}^{p+1}_{L^{p+1}}
\\
&=(\l^2-\l^{\frac{p-1}{2}})\frac{p-1}{2(p+1)}\nr{v}^{p+1}_{L^{p+1}},
\end{align*}
which is positive for all $\l\in(0,1)$, since $p>5$. Hence, {\it (2)} follows. {\it (3)} follows form  simple calculation:
\begin{align*}
\frac{\partial}{\partial \l} S(v_\l)&=\l\nr{v'}^2_{L^2}-\l c\nr{\frac{v}{x}}^2_{L^2}-\l^{\frac{p-1}{2}-1}\frac{p-1}{2(p+1)}\nr{v}^{p+1}_{L^{p+1}}
\\
&=\frac{1}{\l}Q(v_\l).
\end{align*}
To show {\it (4)}, we note that
\[
Q(v_\l)=\frac{\l^2}{(\l^*)^2} Q(v_{\l^*}) + \l^2\left( (\l^*)^{\frac{p-5}{2}} -\l^{\frac{p-5}{2}} \right)\frac{p-1}{2(p+1)}\nr{v}^{p+1}_{L^{p+1}}.
\]
Since $p>5$ and $Q(v_{\l^*})=0$, we get that $\l>\l^*$ implies $Q(v_\l)<0$, and $\l<\l^*$ implies $Q(v_\l)>0$. This and {\it (3)}, implies {\it (4)}.

Finally, we get by simple calculation that
\[
\frac{\partial^2}{\partial \l^2} S(v_\l)= \frac{1}{\l^2} Q(v_\l)-\l^{\frac{p-5}{2}}\left(\frac{p-1}{2}-2 \right)\frac{p-1}{2(p+1)}\nr{v}^{p+1}_{L^{p+1}}.
\]
Since $p>5$, we obtain for $\l>\l^*$ that $\frac{\partial^2}{\partial \l^2} S(v_\l)<0$ which concludes the proof of {\it (5)}.
\end{proof}

To prove orbital instability we prove a new variational characterization of the ground state. Let us define the following set
\[
\mathcal{M}=\{v\in H^1_0(\R^+)\setminus \{0\}: Q(v)=0, J(v)\leq 0\},
\]
and the corresponding minimization problem
\[
d=\inf_{W\in \mathcal{M}} S(W).
\]

Then we have the following.

\begin{lemma} \label{characterization}
The following equality holds:
\[
m=d,
\]
where $m$ is defined by \eqref{minimization1}.
\end{lemma}

\begin{proof}
 Let $v\in \mathcal{G}$. Since $v$ solves \eqref{statIV}, by Lemma \ref{identity} we have that $Q(v)=J(v)=0$, hence $\mathcal{G}\subset \mathcal{M}$, and
\[
d \leq m.
\]
Let now $v\in \mathcal{M}$. Assume first, that $J(v)=0$. In this case $v\in \mathcal{N}$, and $m\leq S(v)$. Let us assume that $J(v)<0$. Then for $v_\l(x)=\l^{1/2}v(\l x)$ we have
\[
J(v_\l)=\l^2\nr{v'}^2_{L^2}-\l^2c\nr{\frac{v}{x}}^2_{L^2}+\o\nr{v}^2_{L^2}-\l^{(p-1)/2}\nr{v}^{p+1}_{L^{p+1}},
\]
and $\lim_{\l\downarrow 0} J(v_\l)>\o \nr{v}^2_{L^2}$, thus there exists $\l_1\in (0,1)$, such that $J(v_{\l_1})=0$. By Proposition \ref{existenceIV} 
\[
m\leq S(v_{\l_1}).
\]
From $Q(v)=0$ and Lemma \ref{scaling} we have
\[
S(v_{\l_1})\leq S(v),
\]
hence $m\leq S(v)$ for all $v\in \mathcal{M}$. Therefore $m\leq d$, which concludes the proof.
\end{proof}

We now define the manifold
\[
\mathcal{J}=\{u\in H^1_0(\R^+)\setminus \{0\}: J(u)<0, Q(u)<0, S(u)< d\}.
\]
We will prove the invariance of $\mathcal{J}$ under the flow of \eqref{eq-wave}.

\begin{lemma} \label{invariance}
Let $u_0 \in \mathcal{J}$ and $u\in C((-T_{\mathrm{min}},T_{\mathrm{max}}), H^1_0(\R^+))$ the corresponding solution to \eqref{eq-wave}. Then $u(t)\in \mathcal{J}$ for all $t\in (-T_{\mathrm{min}}, T_{\mathrm{max}})$.
\end{lemma}

\begin{proof}
Let $u_0\in \mathcal{J}$ and $u\in C((-T_{\mathrm{min}}, T_{\mathrm{max}}), H^1_0(\R^+))$ the corresponding maximal solution. Since $S$ is conserved under the flow of \eqref{eq-wave} we have for all $t\in (-T_{\mathrm{min}}, T_{\mathrm{max}})$ that
\[
S(u(t))=S(u_0)<d.
\]
We prove the assertion by contradiction. Suppose that there exists $t\in (-T_{\mathrm{min}}, T_{\mathrm{max}})$ such that
\[
J(u(t))\geq 0.
\]
Then, since $J$ and $u$ are continuous, there exists $t_0\in (-T_{\mathrm{min}}, T_{\mathrm{max}})$ such that
\[
J(u(t_0))=0,
\]
thus $u(t_0) \in \mathcal{N}$. Then by Proposition \ref{existenceIV} we have that
\[
S(u(t_0))\geq d,
\]
which is a contradiction, thus $J(u(t))<0$ for all $t \in (-T_{\mathrm{min}},T_{\mathrm{max}})$. Let us suppose now that for some $t\in (-T_{\mathrm{min}}, T_{\mathrm{max}})$ we have
\[
Q(u(t))\geq 0.
\]
Again, by continuity, there exists $t_1 \in(-T_{\mathrm{min}}, T_{\mathrm{max}})$ such that
\[
Q(u(t_1))=0.
\]
Hence we that $Q(u(t_1))=0$, and $J(u(t_1))<0$. Therefore, by Lemma \ref{characterization}
\[
S(u(t_1))\geq d,
\]
which is a contradiction. Hence,
\[
Q(u(t))<0 
\]
for all $t \in (-T_{\mathrm{min}},T_{\mathrm{max}})$, which concludes the proof. 
\end{proof}

\begin{lemma} \label{negative}
Let $u_0 \in \mathcal{J}$ and $u\in C((-T_{\mathrm{min}},T_{\mathrm{max}}),H^1_0(\R^+))$. Then there exists $\eps>0$ such that $Q(u(t))\leq -\eps$ for all $t\in (-T_{\mathrm{min}},T_{\mathrm{max}})$. 
\end{lemma}

\begin{proof}
Let $u_0 \in \mathcal{J}$ and let us define $v:=u(t)$ and $v_\l(x)=\l^{1/2}v(\l x)$. By Lemma \ref{scaling}, there exists $\l_0<1$ such that $Q(v_{\l^*})=0$. If $J(v_{\l^*})\leq 0$, then by Lemma \ref{invariance} we get $S(v_{\l^*})\geq m.$ On the other hand, if $J(v_{\l^*})>0$, there exists $\l_1\in (\l^*,1)$, such that $J(\l_1)=0$ and we replace $\l^*$ with $\l_1$. In this case, by Lemma \ref{characterization} we get $S(v_{\l^*})\geq m$. In conclusion, in both cases we obtain
\begin{equation} \label{neg1}
S(v_{\l^*})\geq d.
\end{equation}
By Lemma \ref{scaling} we know that $\l \mapsto S(v_\l)$ is concave on $(\l^*, +\infty)$, thus
\begin{equation}
S(v)-S(v_{\l_0})\geq (1-\l_0)\frac{\partial}{\partial \l} S(v_\l)\Big|_{\l=1}.
\end{equation}
From Lemma \ref{scaling} we have
\begin{equation}
\frac{\partial}{\partial \l}S(v_\l)\Big|_{\l=1}=Q(v).
\end{equation}
Moreover, since $Q(v)<0$ and $\l^*\in (0,1)$, we have
\begin{equation} \label{neg2}
(1-\l^*)Q(v) > Q(v).
\end{equation}
Combining \eqref{neg1}-\eqref{neg2}, we obtain
\[
S(v)-d > Q(v).
\]
Define $-\eps= S(v)-d$. Then $\eps>0$, since $v\in \mathcal{J}$. Owing to the conservation of the energy and mass, $\eps>0$ is independent from $t$, which concludes the proof.
\end{proof}

\begin{lemma} \label{blow-upIV}
Let us take $u_0\in \mathcal{J}$ such that $x u_0 \in L^2(\R^+)$. Then the maximal solution $u\in C((-T_{\mathrm{min}},T_{\mathrm{max}}),H^1_0(\R^+))$ corresponding to the initial value problem \eqref{eq-wave} blows up in finite time.
\end{lemma}

\begin{proof}
From Lemma \ref{negative} we know that there exists $\eps>0$ such that 
\[
Q(u(t))< -\eps \text{ for } t\in (-T_{\mathrm{min}},T_{\mathrm{max}}).
\]
From Proposition \ref{virial} we know that $\frac{\partial^2}{\partial t^2}\nr{xu(t)}^2_{L^2}=8Q(u(t))$, and by integration we get
\begin{equation} \label{blowup1}
\nr{xu(t)}^2_{L^2}\leq -4\eps t^2 + C_1t+C_2.
\end{equation}
The right hand side of \eqref{blowup1} is negative for large $|t|$, which contradicts with $\nr{x u(t)}^2_{L^2}>0$ for all $t$. Therefore, $T_{\mathrm{min}}>-\infty$ and $T_{\mathrm{max}}<\infty$ and by local well-posedness it follows that
\[
\lim_{t\downarrow -T_{\mathrm{min}}}\nr{u(t)}_{H^1}=+\infty, \text{ and } \lim_{t\uparrow T_{\mathrm{max}}} \nr{u(t)}_{H^1}=+\infty.
\]
\end{proof}

\begin{proof}[Proof of Proposition \ref{unstable}] Let $\ffi \in\mathcal{G}$. Owing to Lemma \ref{blow-upIV}, it suffices to show that there exists a sequence $\{\ffi_\l\}\subset\mathcal{J}$, which converges to $\ffi$ in $H^1_0(\R^+)$. Let us put $\ffi_\l(x)=\l^{1/2}\ffi(\l x)$. By Lemma \ref{scaling} $\{\ffi_\l\}\subset\mathcal{J}$ for all $\l \in (0,1)$. Additionally, by Proposition \ref{regularity}, $\ffi$ decays exponentially at infinity, and so does $\ffi_\l$. Therefore, $x\ffi_\l \in L^2(\R^+)$. Clearly, $\ffi_\l\rightarrow \ffi$ as $\l\rightarrow0$, and by Lemma \ref{blow-upIV} the maximal solution of \eqref{eq-wave} corresponding to $\ffi_\l$, blows up in finite time for all $\l\in (0,1)$. Hence, the conclusion follows.
\end{proof}


\section{Appendix}

We prove the following Lemma:

\begin{lemma} \label{approximations}
Let $\psi_A(x)=q(x+A)-q(x-A)$, where $q$ is \eqref{explicitIV}. Then $\psi_A\in H^1_0(\R^+)$ and for large $A>0$, we have the following approximations:
\begin{align}
\int_0^\infty |\psi_A'|^2dx&= \int_{-\infty}^\infty|q'|^2dx+O\left(\left(2A+\frac{1}{\sqrt{\o}}\right)e^{-2\sqrt{\o} A}\right),\label{appr1}
\\ 
\int_0^\infty |\psi_A|^2dx&= \int_{-\infty}^\infty |q|^2dx +O\left(\left(2A+\frac{1}{\sqrt{\o}}\right)e^{-2\sqrt{\o} A}\right),\label{appr2}
\\ 
\int_0^\infty\frac{|\psi_A(x)|^2}{x^2}&\lesssim\frac{1}{A^2}\int_{-\infty}^\infty |q|^2dx+ O\left(\frac{1}{A^2}e^{-\sqrt{\o} A}\right), \label{appr3}
\\ 
\int_0^\infty|\psi_A(x)|^{p+1}dx &=\int_{-\infty}^\infty |q|^{p+1}dx + O(e^{-2\sqrt{\o} A}). \label{appr4}
\end{align}
\end{lemma}

\begin{proof}
We will use the fact that $q(x)\leq Me^{-\sqrt{\o}|x|}$ and $q'(x)\leq Me^{-\sqrt{\o}|x|}$ for some $M>0$. 

We get \eqref{appr1} by using the symmetry of $q$ and $q'$:
\[
\int_0^\infty |\psi_A'|^2dx=\int_{-\infty}^\infty |q'|^2dx -\int_{-\infty}^\infty q'(x+A)q'(x-A)dx.
\]
We estimate the second term by
\[
\left|\int_{-\infty}^\infty q'(x+A)q'(x-A)dx\right|\lesssim \int_{-\infty}^\infty e^{-\sqrt{\o}|x+A|-\sqrt{\o}|x-A|}dx=\left(\left(2A+\frac{1}{\sqrt{\o}}\right)e^{-2\sqrt{\o} A} \right),
\]
hence \eqref{appr1} follows. We get \eqref{appr2} the same way.

We now show \eqref{appr3}. From Hardy's inequality we get
\begin{align*}
\int_0^{A/2}\frac{|\psi_A|^2}{x^2}dx\leq 4\int_0^{A/2}|\psi'_A(x)|^2=O(e^{-\sqrt{\o} A}).
\end{align*}
Moreover, we have
\begin{align*}
\int_{A/2}^\infty \frac{|\psi_A(x)|^2}{x^2}dx \leq \frac{4}{A^2}\int_{A/2}^\infty |\psi_A|^2dx=\frac{4}{A^2}\int_{-\infty}^\infty |q|^2+O\left( \frac{1}{A^2}e^{-\sqrt{\o} A}\right).
\end{align*}
Hence
\begin{align*}
\int_0^\infty \frac{|\psi_A|^2}{x^2}dx=\int_0^{A/2}\frac{|\psi_A|^2}{x^2}dx+\int^\infty_{A/2}\frac{|\psi_A|^2}{x^2}dx \leq \frac{4}{A^2}\int_{-\infty}^\infty|q|^2dx+O\left(\frac{1}{A^2}e^{-\sqrt{\o} A}\right),
\end{align*}
which is the estimate in \eqref{appr3}.

To show \eqref{appr4}, we use the fact that 
\begin{align*}
|q(x-A)-q(x+A)|^{p+1}=q^{p+1}(x-A)-(p+1)q^p(x-A)q(x+A)+O(q^2(x+A)).
\end{align*}
We get
\begin{align*}
&\int_0^\infty q^{p+1}(x-A)dx =\int_{-\infty}^\infty q^{p+1}(x)dx-\int_{-\infty}^{-A}q^{p+1}(x)dx=\int_{-\infty}^\infty q^{p+1}(x)dx+O(e^{-\sqrt{\o}(p+1)A}),
\\
&\int_0^\infty q^p(x-A)q(x+A)dx\lesssim\int_0^\infty e^{-\sqrt{\o} p|x-A|-\sqrt{\o}|x+A|}dx
=O\left(e^{-2\sqrt{\o} A}\right),
\\
&\int_0^\infty O(q^2(x+A))dx = O(e^{-2\sqrt{\o} A}).
\end{align*}
Hence
\begin{align*}
\int_0^\infty |\psi_A(x)|^{p+1}dx= \int_{-\infty}^\infty |q|^{p+1}dx+O(e^{-2\sqrt{\o} A}).
\end{align*}
This concludes the proof.
\end{proof}

We now state the proof of Lemma \ref{decomposition}. The proof follows the arguments of the paper \cite{JeanjeanIV}, with some important modifications. We introduce the norm
\[
\nr{u}^2=\int_0^\infty \left(|u'|^2-c\frac{|u|^2}{x^2}+\om |u|^2\right)dx,
\]
which is equivalent to the standard norm on $H^1_0(\R^+)$ if $0<c<1/4$.
\begin{proof}{\it Proof of Lemma \ref{decomposition}}

{\it Step 1. There exists $u_0\in H^1_0(\R^+)$, such that, up to a subsequence, $u_n$ is weakly convergent to $u_0$ in $H^1_0(\R^+)$, and $S'(u_0)=0.$}
\\
Since $\{u_n\}_{n\in\N}$ is bounded in $H^1_0(\R^+)$, it admits a weakly convergent subsequence in $H^1_0(\R^+)$ with a weak limit $u_0\in H^1_0(\R^+)$. We only need to show that $S'(u_0)=0$. Since by our assumption $S'(u_n)\rightarrow 0$, it suffices to show that for all $\ffi \in C^\infty_0(\R^+)$ we have
\[
S'(u_n)\ffi-S'(u_0)\ffi \rightarrow 0.
\]
Indeed, we have
\begin{align*}
S'(u_n)\ffi-S'(u_0)\ffi= &\Re\int^\infty_0 (u_n'-u_0')\bar{\ffi}'dx-c\Re\int_0^\infty\frac{(u_n-u_0)\bar{\ffi}}{x^2}dx+\om\Re\int_0^\infty(u_n-u_0)\bar{\ffi} dx\\
&-\Re\int_0^\infty(|u_n|^{p-1}u_n-|u_0|^{p-1}u_0)\bar{\ffi} dx.
\end{align*}
Since $u_n\rightharpoonup u_0$ in $H^1_0(\R^+)$ and strongly in $L^q_{\mathrm{loc}}(\R^+)$ for all $q\geq 1$, our statement follows.
\\
\\
Let us set $v_n=u_n-u_0$.
\\
\\
{\it Step 2. Assume that
\begin{equation}\label{condition1}
\sup_{z\in \R^+}\int_{B_1(z)}|v_n|^2dx\rightarrow 0,
\end{equation}
where $B_1(z)$ is the unit ball centered at $z$.
Then $u_n\rightarrow u_0$ strongly in $H^1_0(\R^+)$, and Lemma \ref{decomposition} holds with $k=0$.}
\\
Using the fact that $S'(u_0)=0$, we get
\begin{align*}
S'(u_n)v_n&=\Re\int_0^\infty u'_n\bar{v}'_n dx -c\Re\int^\infty_0 \frac{u_n\bar{v}_n}{x^2}dx+\om\Re\int_0^\infty u_n\bar{v}_ndx-\Re\int_0^\infty|u_n|^{p-1}u_n\bar{v}_ndx=\\
&=\nr{v_n}^2+\Re\int_0^\infty (|u_0|^{p-1}u_0-|u_n|^{p-1}u_n)\bar{v}_n dx.
\end{align*}
Hence,
\begin{align*}
\nr{v_n}^2=S'(u_n)v_n+\Re\int^\infty_0(|u_n|^{p-1}u_n-|u_0|^{p-1}u_0)\bar{v}_ndx.
\end{align*}
We recall that $S'(u_n)\rightarrow 0$. H\"older's inequality implies that
\[
\left|\int_0^\infty |u_n|^{p-1}u_nv_n dx\right|\leq \nr{u_n}_{L^{p+1}}^p\nr{v_n}_{L^{p+1}}.
\]
Assumption \eqref{condition1} and Lemma 1.1 in \cite{Lions2IV} implies that $\nr{v_n}_{L^{p+1}}\rightarrow 0$. Hence
\[
\Re\int_0^\infty |u_n|^{p-1}u_n\bar{v}_n dx \rightarrow 0.
\]
We obtain similarly that $\Re\int_0^\infty |u_0|^{p-1}u_0\bar{v}_n dx \rightarrow 0$, hence $\nr{v_n}^2\rightarrow 0$, which completes the proof of Step 2.
\\
\\
{\it Step 3. Assume that there exist $\{z_n\}_{n\in\N}\subset \R^+$ and $d>0$, such that
\begin{equation}\label{condition2}
\int_{B_1(z_n)}|v_n|^2 dx \rightarrow d.
\end{equation}
Then, up to a subsequence, we have for $q\in H^1(\R)$, that
(i) $z_n\rightarrow \infty$,
(ii) $u_n(\cdot+z_n)\rightharpoonup q \neq 0$ in $H^1(\R)$, and
(iii) ${S^\infty}'(q)=0$.}

To show {\it(i)}, let us assume by contradiction that $\{z_n\}_{n\in\N}$ has an accumulation point $z^*\in \R^+$. Then for a subsequence of $\{v_n\}_{n\in\N}$ we have
\[
\int_{B_2(z^*)}|v_n|^2dx \geq d.
\]
Since $v_n\rightharpoonup 0$ in $H^1_0(\R^+)$, we have $v_n \rightarrow 0$ in $L^2(B_2(z^*))$, which implies that
\[
d \leq \lim_{n\rightarrow \infty}\int_{B_2(z^*)}|v_n|^2dx=0,
\]
which is a contradiction, hence {\it (i)} holds.

Since $u_n(\cdot + z_n)$ is bounded in $H^1(\R)$ the re exists $q\in H^1(\R)$ such that $u_n(\cdot+z_n)$ converges weakly to $q$ in $H^1(\R)$. We only need to show that $q\neq 0$. Since $u_0(\cdot + z_n) \rightharpoonup 0$ in $H^1(\R)$, we have that $v_n(\cdot + z_n) \rightharpoonup q$ in $H^1(\R)$, and in $L^2_{\mathrm{loc}}(\R)$ in particular. Hence
\[
\int_{B_1(0)}|q(x)|^2dx=\lim_{n\rightarrow \infty}\int_{B_1(0)}|v_n(x+z_n)|^2dx=\int_{B_1(z_n)}|v_n(y)|^2dy\geq d>0.
\]
This implies that $q\neq 0$.

We finally show {\it(iii)}. We define $\tilde{u}(\cdot)=u_n(\cdot +z_n)$. We obtain, similarly as in Step 1, that for any $\ffi\in C^\infty_0(\R)$,
\[
{S^\infty}'(\tilde{u}_n)\ffi-{S^\infty}'(q)\ffi\rightarrow 0.
\]
It remains to show that ${S^\infty}'(\tilde{u}_n)\ffi\rightarrow 0$. For any fixed $\ffi \in C^\infty_0(\R)$, $\ffi(\cdot-z_n)$ is in $H^1_0(\R^+)$ for sufficiently big $n\in\N$. Hence, we obtain
\begin{align*}
S'(u_n)\ffi(\cdot-z_n)&=\Re\int_{-z_n}^\infty u'_n(x+z_n)\bar{\ffi}'_n(x) dx-c\Re\int_{-z_n}^\infty\frac{u_n(x+z_n)\bar{\ffi}(x)}{(x+z_n)^2}dx\\
&+\om\Re\int_{-z_n}^\infty u_n(x+z_n)\bar{\ffi}(x)dx-\Re\int_{-z_n}^\infty |u_n(x+z_n)|^{p-1}u_n(x+z_n)\bar{\ffi}(x)dx. 
\end{align*}
Since $S'(u_n)\rightarrow 0$ and $\ffi(\cdot - z_n)$ is bounded in $H^1(\R)$, it follows
\begin{align*}
\Re\int_{-z_n}^\infty &\tilde{u}'_n(x)\bar{\ffi}'_n(x) dx-c\Re\int_{-z_n}^\infty\frac{\tilde{u}_n(x)\bar{\ffi}(x)}{(x+z_n)^2}dx
\\
&+\om\Re\int_{-z_n}^\infty \tilde{u}_n(x)\bar{\ffi}(x)dx-\Re\int_{-z_n}^\infty |\tilde{u}_n(x)|^{p-1}\tilde{u}_n(x)\bar{\ffi}(x)dx\rightarrow 0.
\end{align*}
Moreover, since $u_n$ is bounded in $L^\infty$, and $\ffi$ is compactly supported, we get
\begin{align*}
\left|\Re\int_{-z_n}^\infty\frac{\tilde{u}_n(x)\bar{\ffi}(x)}{(x+z_n)^2}dx\right| = \left|\Re\int_0^\infty\frac{u_n(x)\bar{\ffi}(x-z_n)}{x^2}dx\right|\leq \frac{1}{(z_n-\inf\{\supp(\ffi)\})^2}\nr{u_n\ffi}_{L^\infty}\rightarrow 0,
\end{align*}
Thus
\begin{align*}
{S^\infty}'(\tilde{u}_n)\ffi&=\Re\int_{-\infty}^\infty \tilde{u}'_n(x)\bar{\ffi}'_n(x) dx+\om\Re\int_{-\infty}^\infty \tilde{u}_n(x)\bar{\ffi}(x)dx
\\
&-\Re\int_{-\infty}^\infty |\tilde{u}_n(x)|^{p-1}\tilde{u}_n(x)\bar{\ffi}(x)dx\rightarrow 0,
\end{align*}
which concludes the proof of Step 3.
\\
\\
{\it Step 4. Suppose there exist $k\geq 1$, $\{ x_n^i\}\subset \R^+$, $q_i\in H^1(\R)$ for $1\leq i\leq k$, such that
\begin{align*}
x_n^i\rightarrow \infty, \quad |x_n^i-x_n^{j}|\rightarrow \infty \textit{ if } i\neq j,\\
u_n(\cdot + x_n^i)\rightarrow q_i\neq 0, \textit{ for all } 1\leq i \leq k,\\
{S^\infty}'(q_i)=0.
\end{align*}
Then
\\
1) If $\sup_{z\in\R^+}\int_{B_1(z)}|u_n-u_0-\sum_{i=1}^kq_i(\cdot -x_n^i)|^2dx\rightarrow 0$ then
\[
\nr{u_n-u_0-\sum_{i=1}^kq_i(\cdot -x_n^i)}_{H^1}\rightarrow 0.
\]
2) If there exist $\{z_n\}\subset\R^+$ and $d>0$, such that
\[
\int_{B_1(z_n)}\left|u_n-u_0-\sum_{i=1}^kq_i(\cdot -x_n^i)\right|^2dx \rightarrow d,
\]
then, up to a subsequence, it follows that
\begin{align*}
(i)\textit{  } z_n\rightarrow \infty, \textit{ and } |z_n-x_n^i|\rightarrow \infty \textit{ for all } 1\leq i \leq k,\\
(ii)\textit{  } u_n(\cdot +z_n )\rightharpoonup q_{i+1} \quad  (iii)\textit{  } {S^\infty}'(q_{i+1})=0.
\end{align*}
}

Suppose assumption {\it 1)} holds. We introduce $\xi_n=u_n-u_0-\sum_{i=1}^kq^a_i(\cdot - x_n^i)$, where $q^a_i$ is a suitable cut-off of $q_i$, such that $\supp(q^a_i)\subset (0,\infty)$. This is possible owing to the exponential decay of $q_i$ at infinity, and $x_n^i\rightarrow \infty$ as $n\rightarrow \infty$ for all $i$. We get
\begin{align*}
S'(u_n)\xi_n&=\Re\int_0^\infty u'_n\bar{\xi}'_ndx-c\Re\int_0^\infty\frac{u_n\bar{\xi}_n}{x^2}dx+\om\Re\int_0^\infty u_n\bar{\xi}_ndx-\Re\int_0^\infty|u_n|^{p-1}u_n\bar{\xi}_ndx\\
&=\nr{\xi_n}^2+\Re\int_0^\infty(u'_0+\sum_{i=1}^k {q^a_i}'(\cdot-x_n^i))\bar{\xi}'_ndx
\\
 &+\Re\int_0^\infty\left(\om-\frac{c}{x^2}\right)\left(u_0+\sum_{i=1}^kq^a_i(\cdot-x_n^i)\right) \bar{\xi}_ndx-\Re\int_0^\infty|u_n|^{p-1}u_n\bar{\xi}_ndx. 
\end{align*}
Since $S'(u_0)\xi_n=0$, we get
\begin{align*}
S'(u_n)\xi_n&=\nr{\xi_n}^2+\Re\int_0^\infty (|u_0|^{p-1}u_0-|u_n|^{p-1}u_n)\bar{\xi}_ndx+\Re\int_0^\infty\sum_{i=1}^k{q^a_i}'(\cdot-x_n^i)\bar{\xi}'_ndx
\\
&+\Re\int_0^\infty \left(\om-\frac{c}{x^2} \right)\sum_{i=1}^kq^a_i(\cdot-x_n^i)\bar{\xi}_ndx.
\end{align*}
Using the fact that $\nr{\xi_n}_{L^{p+1}}\rightarrow 0$ by Lemma 1.1 in \cite{Lions2IV}, we get that the second term of the right hand side converges to zero. Now, from the weak convergence of $\xi_n$ to zero and that $S'(u_n)\rightarrow 0$,we obtain that $\nr{\xi_n}\rightarrow 0$.

Suppose now that assumption {\it 2)} holds. Then {\it (i)} and {\it (ii)} follows as in Step 3. To show {\it (ii)}, let us set $\tilde{u}_n=u_n(\cdot + z_n)$. We note that
\begin{align*}
{S^\infty}'(\tilde{u}_n)\ffi - {S^\infty}'(q)\ffi\rightarrow 0,
\end{align*}
for all $\ffi\in C^\infty_0(\R)$. Now ${S^\infty}'(\tilde{u}_n)\rightarrow 0$ follows similarly as in Step 3, which concludes the proof.
\\
\\
{\it Step 5. Conclusion}
By Step 1 we know that $u_n\rightharpoonup u_0$ and $S'(u_0)=0$. Hence (i) of Lemma \ref{decomposition} is verified. If the assumption of Step 2 holds, then Lemma \ref{decomposition} is true with $k=0$. Otherwise, the assumption of Step 3 holds. 
We have to iterate Step 4.  We only need to show that assumption 1 of Step 4 occurs after a finite number of iterations. 
Let us notice that
\[
\nr{u_n-u_0-\sum_{i=1}^kq_i(\cdot- x_i^n)}^2_{H^1}=\nr{u_n}^2_{H^1}+\nr{u_0}^2_{H^1}+\sum_{i=1}^k\nr{q_i}^2_{H^1}-2\left\langle u_n,u_0+\sum_{i=1}^kq_i(\cdot-x_i^n)\right\rangle_{H^1}.
\]
Moreover, since $u_n\rightharpoonup u_0$ and $u_n(\cdot + x_i^n)\rightharpoonup q_i$, we get for the last term that
\[
\left\langle u_n,u_0+\sum_{i=1}^kq_i(\cdot-x_i^n)\right\rangle_{H^1} \rightarrow \nr{u_0}^2_{H^1}+\sum_{i=1}^k\nr{q_i}^2_{H^1},
\]

Now since $u_n$ converges weakly to $u_0$, we obtain for $k\geq 1$ that
\[
\lim_{n\rightarrow \infty}\nr{u_n}^2_{H^1}-\nr{u_0}^2_{H^1}-\sum_{i=1}^k\nr{q_i}^2_{H^1}=\lim_{n\rightarrow \infty}\nr{u_n-u_0-\sum_{i=1}^kq_i(\cdot - x_i^n)}^2_{H^1}\geq 0.
\]
Since $q_i$ is a nontrivial critical point of $S^\infty$, it is true that $\nr{q_i}_{H^1}\geq \epsilon>0$. Hence, after a finite number of iterations assumption 1 of Step 4 must occur.

Finally, we have to verify that
\[
S(u_n)\rightarrow S(u_0)+\sum_{i=1}^kS^\infty(q_i).
\]
We first show that
\begin{equation} \label{first}
S(u_n)\rightarrow S(u_0)+S^\infty(v_n).
\end{equation}
A straightforward calculation gives
\begin{align*}
S(u_n)&=S(u_0)+S^\infty(v_n)+\Re\int_0^\infty u'_0(\bar{u}'_n-\bar{u}'_0)dx-c\Re\int_0^\infty\frac{u_0(\bar{u}_n-\bar{u}_0)}{x^2}dx
\\
&+\o\Re\int_0^\infty u_0(\bar{u}_n-\bar{u}_0)dx
-\frac{c}{2}\int_0^\infty\frac{|u_n-u_0|^2}{x^2}dx
\\
&+\frac{1}{p+1}\left( \nr{u_n-u_0}^{p+1}_{L^{p+1}}-\nr{u_n}^{p+1}_{L^{p+1}}+\nr{u_n}^{p+1}_{L^{p+1}}\right)
\end{align*}
From a lemma by Brezis and Lieb (see e.g. Lemme 4.6 \cite{KavianIV}) we have
\[
\int^\infty_0|u_n-u_0|^{p+1}dx-\int^\infty_0|u_n|^{p+1}dx+\int^\infty_0|u_0|^{p+1}dx\rightarrow 0.
\]
Hence \eqref{first} follows. It only remains to show that
\[
S^\infty(v_n)\rightarrow \sum_{i=1}^kS^\infty(q_i).
\]
We calculate
\begin{align*}
S(v_n)&=\frac{1}{2}\nr{v_n-\sum_{i=1}^kq_i(\cdot - x_i^n)}^2_{H^1}+\frac{1}{2}\nr{\sum_{i=1}^kq_i(\cdot - x_i^n)}^2_{H^1}
\\
&+\left\langle v_n-\sum_{i=1}^kq_i(\cdot -x_i^n),\sum_{i=1}^kq_i(\cdot - x_i^n) \right\rangle_{H^1}
-\frac{1}{p+1}\nr{\sum_{i=1}^kq_i(\cdot - x_i^n)}^{p+1}_{L^{p+1}}
\\
&-\frac{1}{p+1}\nr{v_n}^{p+1}_{L^{p+1}}+\frac{1}{p+1}\nr{\sum_{i=1}^kq_i(\cdot-x_i^n)}^{p+1}_{L^{p+1}}.
\end{align*}
We have shown that $v_n-\sum_{i=1}^kq_i(\cdot - x_i^n)\rightarrow 0$ strongly in $H^1$. Hence the first and third term above converges to zero as $n\rightarrow\infty$. By using Sobolev's inequality and $\nr{A-B}\geq |\nr{A}-\nr{B}|$ we have
\[
\nr{\sum_{i=1}^kq_i(\cdot-x_i^n)}^{p+1}_{L^{p+1}}-\nr{v_n}^{p+1}_{L^{p+1}}\rightarrow 0,
\]
which concludes the proof. \end{proof}


\begin{thebibliography}{99}
\small

\bibitem{BensouilahIV} A.~Bensouilah, V.~D.~Dinh, S.~Zhu,
On stability and instability of standing waves for the nonlinear Schr\"odinger equation with inverse-square potential,
{\it Journal of Mathematical Physics} 59 (2018).

\bibitem{BruneauIV} L.~Bruneau, J.~Derezi\'nski, V.~Georgescu,
Homogeneous Schr\"odinger operators on half-line,
{\it Annales Henri Poincare} 12(3) (2009), 547--590.

\bibitem{CazenaveIV} T.~Cazenave,
Semilinear Schr\"odinger equations,
{\it American Mathematical Society} (2003).

\bibitem{Cazenave1IV} T.~Cazenave, P.~L.~Lions,
Orbital stability of standing waves for some nonlinear Schr\"odinger equations,
{\it Comm. Math. Phys.} 85 (1982), 594--561. 

\bibitem{CH98IV} T.~Cazenave, A.~Haraux,
An Introduction to Semilinear Evolution Equations,
{\it Oxford Lecture Series in Mathematics and its Applications 13} (1998).


\bibitem{CsoboIV} E.~Csobo, F.~Genoud,
Minimal mass blow-up solutions for the $L^2$ critical NLS with inverse-square potential,
{\it Nonlinear Analysis} 168 (2018), 110--129.

\bibitem{DaviesIV} E.~B.~Davies,
A review of Hardy inequalities,
{\it Operator Theory Advances and Applications} 110(2) (1999), 55--68.

\bibitem{DinhIV} V.~D.~Dinh,
Global existence and blow-up for a class of the focusing nonlinear Schr\"odinger equation with inverse-square potential,
{\it Journal of Mathematical Analysis and Applications} 468(1) (2018), 270--303.



\bibitem{FukuizumiIV} R.~Fukuizumi, L.~Jeanjean,
Stability for standing waves for a nonlinear Schr\"odinger equation with repulsive Dirac delta potentials,
{\it Discrete and Continuous Dynamical Systems} 21(1) (2008), 121--136.


\bibitem{HajaiejIV} H.~Hajaiej, C.~A.~Stuart,
On the variational approach to the stability of standing waves of the nonlinear Schr\"odinger equation,
{\it Adv. Nonlinear Stud.} (2004), 469--501

\bibitem{JeanjeanIV} L.~Jeanjean, K.~Tanaka,
A positive solution for a nonlinear Schr\"odinger equation on $\R^N$,
{\it  Indiana University Mathematics Journal } (2005), 443--464.

\bibitem{KavianIV} O.~Kavian,
Introduction \'a la Th\`eorie des Points Critiques et Applications aux Probl\`emes Elliptiques
{\it Springer Verlag} (1993).

\bibitem{KovarikIV} H.~Kovarik, F.~Truc,
Schr\"odinger operators on a half-line with inverse square potentials,
{\it Mathematical Modeling of Natural Phenomena,} (2014), 170--176.

\bibitem{LeCozIV} S.~Le~Coz,
Standing waves in nonlinear Schr\"odinger equations,
{\it Analytical and Numerical Aspects of Partial Differential Equations} (2009), 151--192.

\bibitem{LionsIV} P.~L.~Lions,
The concentration-compactness principle in the calculus of variations. The limit case, part 1,
{\it Annales de L. H. P., section C} (1984), 109--145.

\bibitem{Lions2IV} P.~L.~Lions,
The concentration-compactness principle in the calculus of variations. The locally compact case, part 2,
{\it Annales de L. H. P., section C} (1984), 223--283.

\bibitem{MontefuscoIV} E.~Montefusco, Lower Semicontinuity of Functionals via the Concentration-Compactness Principle,
{\it J. of Mathematical Analysis and Applications 263,} (2001), 264--276.

\bibitem{OkazawaIV} N.~Okazawa, T.~Suzuki, T.~Yokota,
Energy methods for abstract nonlinear Schr\"odinger equations,
{\it Evol. Equ. Control Theory} 1 (2012),  337--354.

\bibitem{Shioji} N.~Shioji, K.~Watanabe,
A generalized Pohozaev identity and uniqueness of positive 
radial solutions of $\Delta u +g(r)u+h(r)u^p=0$,
{\it J: Differential Equations} 255 (2013),  4448--4475. 

\bibitem{StruweIV} M.~Struwe,
Variational Methods, second edition,
{\it Springer} (1991). 

\bibitem{SuzukiIV} T.~Suzuki,
Nonlinear Schr\"odinger equations with inverse square potentials in two dimensional space,
{\it Dynamical Systems, Differential Equations and Applications,} (2015), 1019--1024.

\bibitem{TrachanasIV} G.~P.~Trachanas, N.~B.~Zographopoulos,
Orbital stability for the Schr\"odinger operator involving inverse square potential, 
{\it J. Differential Equations} 259 (2015), 4989--5016

\bibitem{WeinsteinIV}  M.~I.~Weinstein,
Nonlinear Schr\"odinger equations and sharp interpolation estimates,
{\it Comm. Math. Phys.} 87 (1982), 567--576. 

\end{thebibliography}
\end{document}